# Shape Preserving Zipper Hidden Variable Fractal Interpolation Function

Chol Hui Yun[1]*, Yu Jong Pak[2], Mi Gyong Ri[3], Kyong Ju Ri[4]

[1), 2), 3), 4)] Faculty of Mathematics, **Kim Il Sung** University,

Pyongyang, Democratic People's Republic of Korea

[1]* ch.yun@ryongnamsan.edu.kp

[2] math8@ryongnamsan.edu.kp

[3] mk.ri1004@ryongnamsan.edu.kp

[4] math9@ryongnamsan.edu.kp

## Abstract

In this paper, we study a new class of zipper fractal interpolation functions (ZFIFs) constructed using a zipper hidden variable iterated function system (ZHVIFS). ZFIFs have more diverse shape than usual fractal interpolation functions, and the hidden variable iterated function system is more general than the iterated function system. Firstly, we introduce a ZHVIFS and construct univariate ZFIFs using the ZHVIFS. Next, we find the conditions on vertical scaling factors for the ZFIFs to preserve boundedness, positivity and piecewise slop of the data set, which are demonstrated through examples.

## 1. Introduction

By means of the fractal interpolation functions (FIFs), we can model the objects or phenomena in nature better than smooth functions such as polynomials and spline etc.

In 1986, Barnsley introduced a construction of FIFs on the basis of the theory of the iterated function system (IFS). Since then, diverse constructions of FIFs using the IFS have been presented. Because the objects and phenomena in nature are very complicated and irregular, the most important thing in the theory of fractal interpolation theory is to construct the wider class of diverse FIFs, that is, to find more flexible constructions of FIFs. Therefore, the constructions of FIFs have been developed in two directions.

One is to generalize the IFSs used in the construction from the IFS with constant vertical scaling factors to the hidden variable IFS (HVIFS) with function vertical scaling factors: IFS ([3,4]), recurrent

IFS (RIFS) ([5, 8, 19, 25]), hidden variable IFS (HVIFS) ([6, 7, 10, 11, 22, 26, and 27]), hidden variable RIFS (HVRIFS) ([3, 20, 21]). Since the vertical scaling factors of IFS determine the fractal characteristics of FIFs, they have been developed from constants to functions in all kinds of IFSs. We use the HVIFS with function vertical scaling factors.

The other is to generalize the underlying fixed point theorem from Banach's fixed point theorem to Rakotch's ([13, 15]), Geraghty's ([14, 18]) and Matkowski's one ([2]). We use Banach's fixed point theorem.

In the IFSs mentioned above, homeomorphisms $L_i : I \to I_i$, which map the interval $I \subset \mathbf{R}$ to the $i$th subinterval $I_i$ of the partition of $I$, are all order-preserving. If we allow $L_i$ not only to preserve order but also to reverse order, then we can construct more diverse fractal sets.

On the basis of this idea, in 2003, the authors ([1]) applied the concept of zipper to the construction of fractals, and in [3], zipper fractal curves were constructed and their pointwise regularity was studied. In [23], the researchers studied a construction of zipper fractal Bezier curves and their applications.

Since the theory of fractal interpolation has been widely applied to various fields, many researchers have tried to adapt the methodology of zipper used in constructions of fractal sets to constructing the FIFs ([9, 16, 17, 22, and 24]).

Fractal interpolation surfaces introduced in [16] were special zipper surfaces constructed on the basis of the idea of zipper. In [9], the authors introduced a construction of affine zipper FIFs and studied some properties such as conditions for the graph of zipper FIFs to be contained in a given rectangle. Moreover, in [24], researchers introduced a new technique to approximate a given continuous function by zipper α-fractal functions. The theory of zipper FIF was generalized to the construction of multivariate zipper FIFs ([17, 22]). In constructions of all these zipper FIFs, IFSs were used. We use HVIFS with function vertical scaling factors to construct zipper FIFs.

The shape preserving interpolations have a lot of applications in physics, robotics, computer graphics, computer-aided design, tomography, data visualization, animation, visual space simulation, signal process, control theory, and so on. In most of computer graphics and data visualization, it is essential to construct the interpolation functions whose graphs preserve the boundedness, positivity or piecewise slop of the data and so on.

In [12], the authors studied the sufficient conditions on vertical scaling factors under which the graphs of FIFs have positivity and are laid between two parallel straight lines. In [9], researchers studied the conditions under which the graph of the affine zipper FIF is contained in a prescribed rectangle. Moreover, in [24], the authors investigated an approximation of a given continuous function $f$ by zipper α-FIFs and sufficient conditions on zipper α-FIFs under which these zipper functions preserve the positivity, monotonicity and convexity of the original function $f$. In [22], one kind of zipper fractal interpolation surfaces preserving the positivity was considered. In this paper, we consider the boundedness, positivity and piecewise slop preservation of zipper HVFIFs.

The rest of our paper is organized as follows: In Section 2, we construct zipper hidden variable FIFs using a zipper HVIFS. In Section 3, 4 and 5, we study the conditions on vertical scaling factors under which the constructed zipper hidden variable FIFs are bounded, positive, and preserve piecewise slop, that is, their graphs are restricted by pairs of piecewise lines, respectively.

## 2. Construction of Zipper Hidden Variable Fractal Interpolation Functions

Let us denote the set of natural numbers not greater than $n$ by $\mathbf{N}_n$, and $\mathbf{N}_n \cup \{0\}$ by $\mathbf{N}_n^0$.

First, we recall the concept of zipper iterated function system[9].

**Definition.** Let $X$ be a complete metric space and $\omega_i : X \to X$, $i \in \mathbf{N}_n$ be non-surjective maps. Let vertices $\{u_0, u_1, \cdots, u_n\}$ and a signature $\boldsymbol{\varepsilon} = (\varepsilon_1, \varepsilon_2, \cdots, \varepsilon_n) \in \{0, 1\}^n$ be given. If for any $i \in \mathbf{N}_n$, $\omega_i(u_0) = u_{i-1+\varepsilon_i}$, $\omega_i(u_n) = u_{i-\varepsilon_i}$, then the system $Z = \{X; \omega_i, i \in \mathbf{N}_n\}$ is called a ***zipper iterated function system*** (**ZIFS**) with vertices $\{u_0, u_1, \cdots, u_n\}$ and a signature $\boldsymbol{\varepsilon}$.

If the signature is $\boldsymbol{\varepsilon} = (0, 0, \cdots, 0)$, then ZIFS is the same as the IFS.

Let a data set $P_0$ be given as follows:

$$P_0 = \{(x_i, y_i) \in \mathbf{R}^2 : i \in \mathbf{N}_n^0\} \ (-\infty < x_0 < x_1 < \cdots < x_n < +\infty).$$

We extend $P_0$ to $P = \{(x_i, y_i, z_i) \in \mathbf{R}^3 : i \in \mathbf{N}_n^0\}$, where $z_i$, $i \in \mathbf{N}_n^0$ are parameters. Moreover, we take a signature $\boldsymbol{\varepsilon} = (\varepsilon_1, \varepsilon_2, \cdots, \varepsilon_n) \in \{0, 1\}^n$. We denote $I := [x_0, x_n]$, $I_i := [x_{i-1}, x_i]$, $i \in \mathbf{N}_n$.

We define maps $L_i : I \to I_i$, $i \in \mathbf{N}_n$ by $L_i(x) = a_i x + b_i$, where $a_i = (x_{i-\varepsilon_i} - x_{i-1+\varepsilon_i})/(x_n - x_0)$ and $b_i = (x_{i-1+\varepsilon_i} x_n - x_{i-\varepsilon_i} x_0)/(x_n - x_0)$. Then $L_i$ satisfies $L_i(x_0) = x_{i-1+\varepsilon_i}$ and $L_i(x_n) = x_{i-\varepsilon_i}$.

Let us define maps $\mathbf{F}_i : I \times \mathbf{R}^2 \to \mathbf{R}^2$, $i \in \mathbf{N}_n$ as follows:

$$\mathbf{F}_i(x, y, z) = \begin{pmatrix} F_{i,1}(x, y, z) \\ F_{i,2}(x, y, z) \end{pmatrix}^T = \left[ \begin{pmatrix} p_i(L_i(x)) & q_i(L_i(x)) \\ \tilde{p}_i(L_i(x)) & \tilde{q}_i(L_i(x)) \end{pmatrix} \begin{pmatrix} y \\ z \end{pmatrix} + \begin{pmatrix} r_i(x) \\ \tilde{r}_i(x) \end{pmatrix} \right]^T, \ (x, y, z) \in I \times \mathbf{R}^2,$$

where $p_i, \tilde{p}_i, q_i, \tilde{q}_i : I_i \to \mathbf{R}$, $i \in \mathbf{N}_n$, which are called *vertical scaling factors*, are Lipschitz functions satisfying $\|p_i\|_\infty + \|\tilde{p}_i\|_\infty < 1$ and $\|q_i\|_\infty + \|\tilde{q}_i\|_\infty < 1$, and $r_i, \tilde{r}_i : I \to \mathbf{R}$, $i \in \mathbf{N}_n$ are given as follows:

$$r_i(x) = \left[\frac{x-x_0}{x_n-x_0} y_{i-\varepsilon_i} + \frac{x-x_n}{x_0-x_n} y_{i-1+\varepsilon_i}\right] - p_i(L_i(x))\left[\frac{x-x_0}{x_n-x_0} y_n + \frac{x-x_n}{x_0-x_n} y_0\right]$$

$$- q_i(L_i(x))\left[\frac{x-x_0}{x_n-x_0} z_n + \frac{x-x_n}{x_0-x_n} z_0\right],$$

$$\widetilde{r}_i(x) = \left[\frac{x-x_0}{x_n-x_0} z_{i-\varepsilon_i} + \frac{x-x_n}{x_0-x_n} z_{i-1+\varepsilon_i}\right] - \widetilde{p}_i(L_i(x))\left[\frac{x-x_0}{x_n-x_0} y_n + \frac{x-x_n}{x_0-x_n} y_0\right]$$

$$- \widetilde{q}_i(L_i(x))\left[\frac{x-x_0}{x_n-x_0} z_n + \frac{x-x_n}{x_0-x_n} z_0\right].$$

Then $r_i, \widetilde{r}_i : I \to \mathbf{R}$, $i \in \mathbf{N}_n$ are Lipschitz functions satisfying the following join-up condition:

$$\mathbf{F}_i(x_0, y_0, z_0) = (y_{i-1+\varepsilon_i}, z_{i-1+\varepsilon_i}), \mathbf{F}_i(x_n, y_n, z_n) = (y_{i-\varepsilon_i}, z_{i-\varepsilon_i}).$$

Denote $\mu := \frac{x-x_0}{x_n-x_0}$ $(0 \leq \mu \leq 1)$. Then, we have

$$F_{i,1}(x, y, z) = p_i(L_i(x))\left[y - \left(\frac{x-x_0}{x_n-x_0} y_n + \frac{x-x_n}{x_0-x_n} y_0\right)\right] + q_i(L_i(x))\left[z - \left(\frac{x-x_0}{x_n-x_0} z_n + \frac{x-x_n}{x_0-x_n} z_0\right)\right] +$$

$$+ \left[\frac{x-x_0}{x_n-x_0} y_{i-\varepsilon_i} + \frac{x-x_n}{x_0-x_n} y_{i-1+\varepsilon_i}\right]$$

$$= p_i(L_i(x))\{y - [\mu y_n + (1-\mu)y_0]\} + q_i(L_i(x))\{z - [\mu z_n + (1-\mu)z_0]\} + [\mu y_{i-\varepsilon_i} + (1-\mu)y_{i-1+\varepsilon_i}],$$

$$F_{i,2}(x, y, z) = \widetilde{p}_i(L_i(x))\{y - [\mu y_n + (1-\mu)y_0]\} + \widetilde{q}_i(L_i(x))\{z - [\mu z_n + (1-\mu)z_0]\} + [\mu z_{i-\varepsilon_i} + (1-\mu)z_{i-1+\varepsilon_i}].$$

Now let us denote $E := [-U, U] \times [-\widetilde{U}, \widetilde{U}] \subset \mathbf{R}^2$, where

$$U = \max_{i \in \mathbf{N}_n}\left\{\frac{\|q_i\|_\infty \cdot \|\widetilde{r}_i\|_\infty + (1-\|\widetilde{q}_i\|_\infty) \cdot \|r_i\|_\infty}{(1-\|p_i\|_\infty)(1-\|\widetilde{q}_i\|_\infty) - \|\widetilde{p}_i\|_\infty \cdot \|q_i\|_\infty}\right\}, \quad \widetilde{U} = \max_{i \in \mathbf{N}_n}\left\{\frac{\|\widetilde{p}_i\|_\infty \cdot \|r_i\|_\infty + (1-\|p_i\|_\infty) \cdot \|\widetilde{r}_i\|_\infty}{(1-\|p_i\|_\infty)(1-\|\widetilde{q}_i\|_\infty) - \|\widetilde{p}_i\|_\infty \cdot \|q_i\|_\infty}\right\}.$$

Then $\mathbf{F}_i : I \times E \to E$. In fact, for any $x \in I$, $y \in [-U, U]$ and $z \in [-\widetilde{U}, \widetilde{U}]$, we have

$$|F_{i,1}(x, y, z)| \leq \|p_i\|_\infty \cdot U + \|q_i\|_\infty \cdot \widetilde{U} + \|r_i\|_\infty$$

$$\leq \|p_i\|_\infty \cdot \max_{i \in \mathbf{N}_n}\left\{\frac{\|q_i\|_\infty \cdot \|\widetilde{r}_i\|_\infty + (1-\|\widetilde{q}_i\|_\infty) \cdot \|r_i\|_\infty}{(1-\|p_i\|_\infty)(1-\|\widetilde{q}_i\|_\infty) - \|\widetilde{p}_i\|_\infty \cdot \|q_i\|_\infty}\right\} +$$

$$+\|q_i\|_\infty \cdot \max_{i\in \mathbf{N}_n}\left\{\frac{\|\tilde{p}_i\|_\infty \cdot \|r_i\|_\infty + (1-\|p_i\|_\infty)\cdot \|\tilde{q}_i\|_\infty}{(1-\|p_i\|_\infty)(1-\|\tilde{q}_i\|_\infty)-\|\tilde{p}_i\|_\infty \cdot \|q_i\|_\infty}\right\}+\|r_i\|_\infty$$

$$\leq \max_{i\in \mathbf{N}_n}\left\{\frac{\|p_i\|_\infty \cdot \|q_i\|_\infty \cdot \|\tilde{r}_i\|_\infty + \|p_i\|_\infty \cdot(1-\|\tilde{q}_i\|_\infty)\cdot \|r_i\|_\infty}{(1-\|p_i\|_\infty)(1-\|\tilde{q}_i\|_\infty)-\|\tilde{p}_i\|_\infty \cdot \|q_i\|_\infty}+\right.$$

$$\left.+\frac{\|\tilde{p}_i\|_\infty \cdot \|q_i\|_\infty \cdot \|r_i\|_\infty + (1-\|p_i\|_\infty)\cdot \|q_i\|_\infty \cdot \|\tilde{r}_i\|_\infty}{(1-\|p_i\|_\infty)(1-\|\tilde{q}_i\|_\infty)-\|\tilde{p}_i\|_\infty \cdot \|q_i\|_\infty}+\|r_i\|_\infty\right\}$$

$$=\max_{i\in \mathbf{N}_n}\left\{\frac{\|q_i\|_\infty \cdot \|\tilde{r}_i\|_\infty + (1-\|\tilde{q}_i\|_\infty)\cdot \|r_i\|_\infty}{(1-\|p_i\|_\infty)(1-\|\tilde{q}_i\|_\infty)-\|\tilde{p}_i\|_\infty \cdot \|q_i\|_\infty}\right\}=U$$

and similarly, we get $|F_{i,2}(x, y, z)|\leq \tilde{U}$.

Now, let us define transformations $\omega_i: I\times E\to I_i\times E$, $i=1, 2, \cdots, n$ as follows:

$$\omega_i(x,y,z)=(L_i(x),\ \mathbf{F}_i(x,y,z)),\ x\in I,\ (y,z)\in E.$$

Then $\omega_i$ is well-defined and satisfies

$$\omega_i(x_0,y_0,z_0)=(x_{i-1+\varepsilon_i}, y_{i-1+\varepsilon_i}, z_{i-1+\varepsilon_i}),\ \omega_i(x_n,y_n,z_n)=(x_{i-\varepsilon_i}, y_{i-\varepsilon_i}, z_{i-\varepsilon_i}).$$

and $\mathsf{Z}=\{I\times E;\ \omega_i: i\in \mathbf{N}_n\}$ is called a ***zipper hidden variable iterated function system*** (***ZHVIFS***) with the data set $P$ and the signature $\varepsilon$.

We use the norm $\|g\|_\infty=\max_x |g(x)|$ and denote Lipschitz constant of Lipschitz function $g$ by $l_g$. Moreover, we denote $l_p:=\max_{i\in \mathbf{N}_n}\{l_{p_i}\}$, $l_{\tilde{p}}:=\max_{i\in \mathbf{N}_n}\{l_{\tilde{p}_i}\}$, $l_q:=\max_{i\in \mathbf{N}_n}\{l_{q_i}\}$, $l_r:=\max_{i\in \mathbf{N}_n}\{l_{r_i}\}$, $l_{\tilde{r}}:=\max_{i\in \mathbf{N}_n}\{l_{\tilde{r}_i}\}$, $l_{\tilde{q}}:=\max_{i\in \mathbf{N}_n}\{l_{\tilde{q}_i}\}$, $A:=\max_{i\in \mathbf{N}_n}\{|a_i|\}$, $\bar{S}:=\max_{i\in \mathbf{N}_n}\{\|p_i\|_\infty+\|\tilde{p}_i\|_\infty,\ \|q_i\|_\infty+\|\tilde{q}_i\|_\infty\}<1$.

**Theorem 1**. *If $\bar{S}<1$, then there exists a distance $\rho_\theta$ equivalent to the Euclidean metric such that $\omega_i, i\in \mathbf{N}_n$ are constractive with respect to $\rho_\theta$.*

**Proof.** We take a positive number $\theta$ such that $\theta<\dfrac{1-A}{AU(l_p+l_{\tilde{p}})+A\tilde{U}(l_q+l_{\tilde{q}})+(l_r+l_{\tilde{r}})}$. We introduce the following metric $\rho_\theta: \mathbf{R}^3\times \mathbf{R}^3\to \mathbf{R}$ on $\mathbf{R}^3$:

$$\rho_\theta((x_1, y_1, z_1), (x_2, y_2, z_2)) = |x_1 - x_2| + \theta(|y_1 - y_2| + |z_1 - z_2|), \quad (x_1, y_1, z_1), (x_2, y_2, z_2) \in \mathbf{R}^3.$$

For any $(\hat{x}, \hat{y}, \hat{z})$, $(\tilde{x}, \tilde{y}, \tilde{z}) \in I_i \times E$, we have

$\rho_\theta(\omega_i(\hat{x}, \hat{y}, \hat{z}), \omega_i(\tilde{x}, \tilde{y}, \tilde{z}))$

$= |L_i(\hat{x}) - L_i(\tilde{x})| + \theta(|F_{i,1}(\hat{x}, \hat{y}, \hat{z}) - F_{i,1}(\tilde{x}, \tilde{y}, \tilde{z})| + |F_{i,2}(\hat{x}, \hat{y}, \hat{z}) - F_{i,2}(\tilde{x}, \tilde{y}, \tilde{z})|)$

$= |a_i||\hat{x} - \tilde{x}| + \theta\{|(p_i(L_i(\hat{x}))\hat{y} + q_i(L_i(\hat{x}))\hat{z} + r_i(\hat{x})) - (p_i(L_i(\tilde{x}))\tilde{y} + q_i(L_i(\tilde{x}))\tilde{z} + r_i(\tilde{x}))|$

$\quad + |(\tilde{p}_i(L_i(\hat{x}))\hat{y} + \tilde{q}_i(L_i(\hat{x}))\hat{z} + \tilde{r}_i(\hat{x})) - (\tilde{p}_i(L_i(\tilde{x}))\tilde{y} + \tilde{q}_i(L_i(\tilde{x}))\tilde{z} + \tilde{r}_i(\tilde{x}))|\}$

$\leq A|\hat{x} - \tilde{x}| + \theta|p_i(L_i(\hat{x}))\hat{y} - p_i(L_i(\tilde{x}))\tilde{y}| + \theta|q_i(L_i(\hat{x}))\hat{z} - q_i(L_i(\tilde{x}))\tilde{z}| + \theta|r_i(\hat{x}) - r_i(\tilde{x})|$

$\quad + \theta|\tilde{p}_i(L_i(\hat{x}))\hat{y} - \tilde{p}_i(L_i(\tilde{x}))\tilde{y}| + \theta|\tilde{q}_i(L_i(\hat{x}))\hat{z} - \tilde{q}_i(L_i(\tilde{x}))\tilde{z}| + \theta|\tilde{r}_i(\hat{x}) - \tilde{r}_i(\tilde{x})|.$   (1)

Since $p_i$ is Lipschitz function, we get

$|p_i(L_i(\hat{x}))\hat{y} - p_i(L_i(\tilde{x}))\tilde{y}| \leq |p_i(L_i(\hat{x}))\hat{y} - p_i(L_i(\tilde{x}))\hat{y}| + |p_i(L_i(\tilde{x}))\hat{y} - p_i(L_i(\tilde{x}))\tilde{y}|$

$\quad \leq l_{p_i}|L_i(\hat{x}) - L_i(\tilde{x})| \cdot |\hat{y}| + \|p_i\|_\infty|\hat{y} - \tilde{y}|$

$\quad \leq U l_p A|\hat{x} - \tilde{x}| + \|p_i\|_\infty|\hat{y} - \tilde{y}|$

and similarly, we have

$|q_i(L_i(\hat{x}))\hat{z} - q_i(L_i(\tilde{x}))\tilde{z}| \leq \tilde{U} l_q A|\hat{x} - \tilde{x}| + \|q_i\|_\infty|\hat{z} - \tilde{z}|,$

$|\tilde{p}_i(L_i(\hat{x}))\hat{y} - \tilde{p}_i(L_i(\tilde{x}))\tilde{y}| \leq U l_{\tilde{p}} A|\hat{x} - \tilde{x}| + \|\tilde{p}_i\|_\infty|\hat{y} - \tilde{y}|,$

$|\tilde{q}_i(L_i(\hat{x}))\hat{z} - \tilde{q}_i(L_i(\tilde{x}))\tilde{z}| \leq \tilde{U} l_{\tilde{q}} A|\hat{x} - \tilde{x}| + \|\tilde{q}_i\|_\infty|\hat{z} - \tilde{z}|.$

Therefore, by (1), we obtain

$\rho_\theta(\omega_i(\hat{x}, \hat{y}, \hat{z}), \omega_i(\tilde{x}, \tilde{y}, \tilde{z})) \leq \{A + \theta[AU(l_p + l_{\tilde{p}}) + A\tilde{U}(l_q + l_{\tilde{q}}) + (l_r + l_{\tilde{r}})]\} \cdot |\hat{x} - \tilde{x}| +$

$\quad + \theta[(\|p_i\|_\infty + \|\tilde{p}_i\|_\infty) \cdot |\hat{y} - \tilde{y}| + (\|q_i\|_\infty + \|\tilde{q}_i\|_\infty) \cdot |\hat{z} - \tilde{z}|]$

$\quad \leq \kappa \cdot (|\hat{x} - \tilde{x}| + \theta(|\hat{y} - \tilde{y}| + |\hat{z} - \tilde{z}|))$

$\quad = \kappa \rho_\theta((\hat{x}, \hat{y}, \hat{z}), (\tilde{x}, \tilde{y}, \tilde{z})),$

where $\kappa = \max\{A + \theta(AU(l_p + l_{\tilde{p}}) + A\tilde{U}(l_q + l_{\tilde{q}}) + (l_r + l_{\tilde{r}})), \bar{S}\} < 1$, since $A + \theta(AU(l_p + l_{\tilde{p}}) + A\tilde{M}(l_q + l_{\tilde{q}}) + (l_r + l_{\tilde{r}})) < 1$ and $\bar{S} < 1$. Thus $\omega_i$, $i \in \mathbf{N}_n$ are contractive. □

The set of all non-empty compact subsets of $I \times E$, denoted by $H(I \times E)$, is a complete metric space with respect to Hausdorff metric. Let us define $\mathbf{W}: H(I \times E) \to H(I \times E)$ as follows:

$$\mathbf{W}(B) = \bigcup_{i=1}^{n} \omega_i(B), \quad B \in H(I \times E).$$

If there exists a set $\mathbf{A} \in H(I \times E)$ such that $\mathbf{W}(\mathbf{A}) = \mathbf{A}$ i.e. $\mathbf{A} = \bigcup_{i=1}^{n} \omega_i(\mathbf{A})$, then $\mathbf{A}$ is called an *attractor of the ZHVIFS* $\mathbf{Z}$. By Theorem 1, the ZHVIFS $\mathbf{Z}$ is hyperbolic and hence it has a unique attractor $\mathbf{A}$. Now, we show the existence of a continous map whose graph is $\mathbf{A}$.

**Theorem 2.** *There exists a continous map* $\mathbf{f} = (f_1, f_2)$ *which interpolates the data set* $P$ *and whose graph is* $\mathbf{A}$*, i.e.,* $\mathbf{f}(x_i) = (y_i, z_i)$, $i \in \mathbf{N}_n^0$ *and* $\mathbf{A} = Gr\mathbf{f} = \{(x, \mathbf{f}(x)): x \in I\}$.

**Proof.** The set

$$C_0^2(I) = \{\mathbf{g} = (g_1, g_2) \in C(I) \times C(I) \mid \mathbf{g}(x_0) = (y_0, z_0), \mathbf{g}(x_n) = (y_n, z_n)\}$$

is a complete metric space with respect to the following metric:

$$D(\mathbf{g}, \mathbf{h}) = \max_{x \in I}\{|g_1(x) - h_1(x)| + |g_2(x) - h_2(x)|\}, \quad \mathbf{g} = (g_1, g_2), \mathbf{h} = (h_1, h_2) \in C_0^2(I).$$

Let us define an operator $\mathbf{T}$ on $C_0^2(I)$ as follows: for any $\mathbf{g} \in C_0^2(I)$,

$$(\mathbf{Tg})(x) = \sum_{i=1}^{n} \mathbf{F}_i(L_i^{-1}(x), g_1(L_i^{-1}(x)), g_2(L_i^{-1}(x))) \chi_{I_i}(x), \quad x \in I,$$

where $\chi_{I_i}$ are characteristic functions.

Then, we have $\mathbf{T}: C_0^2(I) \to C_0^2(I)$. In fact, for any $\mathbf{g} \in C_0^2(I)$, we get

$$\lim_{x \to x_i^-}(\mathbf{Tg})(x) = \lim_{x \to x_i^-} \mathbf{F}_i(L_i^{-1}(x), \mathbf{g}(L_i^{-1}(x))) = \mathbf{F}_i(x_{n \cdot (1-\varepsilon_i)}, \mathbf{g}(x_{n \cdot (1-\varepsilon_i)}))$$

$$= \mathbf{F}_i(x_{n \cdot (1-\varepsilon_i)}, y_{n \cdot (1-\varepsilon_i)}, z_{n \cdot (1-\varepsilon_i)}) = (y_i, z_i),$$

$$\lim_{x \to x_i^+}(\mathbf{Tg})(x) = \lim_{x \to x_i^+} \mathbf{F}_{i+1}(L_{i+1}^{-1}(x), \mathbf{g}(L_{i+1}^{-1}(x))) = \mathbf{F}_{i+1}(x_{n \cdot \varepsilon_{i+1}}, \mathbf{g}(x_{n \cdot \varepsilon_{i+1}}))$$

$$= \mathbf{F}_{i+1}(x_{n \cdot \varepsilon_{i+1}}, y_{n \cdot \varepsilon_{i+1}}, z_{n \cdot \varepsilon_{i+1}}) = (y_i, z_i).$$

Hence, we get $\mathbf{Tg} \in C(I) \times C(I)$. Moreover, we have

$$(\mathbf{Tg})(x_0) = \sum_{i=1}^{n} \mathbf{F}_i(L_i^{-1}(x_0), g_1(L_i^{-1}(x_0)), g_2(L_i^{-1}(x_0))) \chi_{I_i}(x_0)$$

$$= \mathbf{F}_1(L_1^{-1}(x_0), g_1(L_1^{-1}(x_0)), g_2(L_1^{-1}(x_0))).$$

Moreover, if $\varepsilon_1 = 0$, then by $L_1^{-1}(x_0) = x_0$, we get

$$(\mathbf{Tg})(x_0) = \mathbf{F}_1(x_0, g_1(x_0), g_2(x_0)) = \mathbf{F}_1(x_0, y_0, z_0) = (y_{\varepsilon_1}, z_{\varepsilon_1}) = (y_0, z_0),$$

and if $\varepsilon_1 = 1$, then by $L_1^{-1}(x_0) = x_n$, we have

$$(\mathbf{Tg})(x_0) = \mathbf{F}_1(x_n, g_1(x_n), g_2(x_n)) = \mathbf{F}_1(x_n, y_n, z_n) = (y_{1-\varepsilon_1}, z_{1-\varepsilon_1}) = (y_0, z_0).$$

Hence, $(\mathbf{Tg})(x_0) = (y_0, z_0)$. Similarly, we can prove $(\mathbf{Tg})(x_n) = (y_n, z_n)$. Therefore, $\mathbf{Tg} \in C_0^2(I)$.

Furthermore, $\mathbf{T}$ is contractive. In fact, for any $\mathbf{g} = (g_1, g_2)$, $\mathbf{h} = (h_1, h_2) \in C_0^2(I)$, we have

$$D(\mathbf{Tg}, \mathbf{Th}) = \max_{x \in I}\{|(\mathbf{Tg})_1(x) - (\mathbf{Th})_1(x)| + |(\mathbf{Tg})_2(x) - (\mathbf{Th})_2(x)|\},$$

and for any $x \in I$, we get

$|(\mathbf{Tg})_1(x) - (\mathbf{Th})_1(x)| =$

$$= \left| \sum_{i=1}^n F_{i,1}(L_i^{-1}(x), g_1(L_i^{-1}(x)), g_2(L_i^{-1}(x))) \chi_{I_i}(x) - \sum_{i=1}^n F_{i,1}(L_i^{-1}(x), h_1(L_i^{-1}(x)), h_2(L_i^{-1}(x))) \chi_{I_i}(x) \right|$$

$$= \sum_{i=1}^n \left| F_{i,1}(L_i^{-1}(x), g_1(L_i^{-1}(x)), g_2(L_i^{-1}(x))) - F_{i,1}(L_i^{-1}(x), h_1(L_i^{-1}(x)), h_2(L_i^{-1}(x))) \right| \chi_{I_i}(x)$$

$$= \sum_{i=1}^n |(p_i(x)g_1(L_i^{-1}(x)) + q_i(x)g_2(L_i^{-1}(x)) + r_i(L_i^{-1}(x))) -$$

$$- (p_i(x)h_1(L_i^{-1}(x)) + q_i(x)h_2(L_i^{-1}(x)) + r_i(L_i^{-1}(x)))| \chi_{I_i}(x)$$

$$\leq \sum_{i=1}^n (\|p_i\|_\infty |g_1(L_i^{-1}(x)) - h_1(L_i^{-1}(x))| + \|q_i\|_\infty |g_2(L_i^{-1}(x)) - h_2(L_i^{-1}(x))|) \chi_{I_i}(x),$$

and

$$|(\mathbf{Tg})_2(x) - (\mathbf{Th})_2(x)| \leq \sum_{i=1}^n (\|\tilde{p}_i\|_\infty |g_1(L_i^{-1}(x)) - h_1(L_i^{-1}(x))| + \|\tilde{q}_i\|_\infty |g_2(L_i^{-1}(x)) - h_2(L_i^{-1}(x))|) \chi_{I_i}(x).$$

Therefore, we have

$|(\mathbf{Tg})_1(x) - (\mathbf{Th})_1(x)| + |(\mathbf{Tg})_2(x) - (\mathbf{Th})_2(x)|$

$$\leq \sum_{i=1}^n [(\|p_i\|_\infty + \|\tilde{p}_i\|_\infty)|g_1(L_i^{-1}(x)) - h_1(L_i^{-1}(x))| + (\|q_i\|_\infty + \|\tilde{q}_i\|_\infty)|g_2(L_i^{-1}(x)) - h_2(L_i^{-1}(x))|] \chi_{I_i}(x)$$

$$\leq \overline{S}\sum_{i=1}^{n}(|g_1(L_i^{-1}(x)) - h_1(L_i^{-1}(x))| + |g_2(L_i^{-1}(x)) - h_2(L_i^{-1}(x))|) \chi_{I_i}(x)$$

and

$$D(\mathbf{Tg}, \mathbf{Th}) \leq \overline{S} D(\mathbf{g}, \mathbf{h}).$$

Then, by Banach fixed point theorem, $\mathbf{T}$ has a unique fixed point $\mathbf{f} = (f_1, f_2)$.

The graph of $\mathbf{f}$ is the attractor of $\mathbf{Z}$. In fact, we get

$$Gr(\mathbf{f}) = \{(x, f_1(x), f_2(x)) \mid x \in I\} = \bigcup_{i=1}^{n}\{(x, \mathbf{F}_i(L_i^{-1}(x), \mathbf{f}(L_i^{-1}(x)))) \mid x \in I_i\}$$

$$= \bigcup_{i=1}^{n}\{(L_i(x), \mathbf{F}_i(x, \mathbf{f}(x))) \mid x \in I\} = \bigcup_{i=1}^{n}\{\omega_i(x, \mathbf{f}(x)) \mid x \in I\}$$

$$= \bigcup_{i=1}^{n} \omega_i(Gr(\mathbf{f})). \quad \square$$

The first component function $f_1 : I \to \mathbf{R}$ of the fixed point $\mathbf{f} = (f_1, f_2)$ of the operator $\mathbf{T}$ interpolates the given data set $P_0$. The interpolation function $f_1$ is called a ***zipper hidden variable fractal interpolation function*** (**ZHVFIF**) and satisfies

$$f_1(x) = \sum_{i=1}^{n}(p_i(x)f_1(L_i^{-1}(x)) + q_i(x)f_2(L_i^{-1}(x)) + r_i(L_i^{-1}(x))) \chi_{I_i}(x), \ x \in I.$$

The second component function $f_2$ is the zipper fractal interpolation function (ZFIF) which interpolates the set $\{(x_i, z_i) \in \mathbf{R}^2 : i \in \mathbf{N}_n^0\}$ and satisfies

$$f_2(x) = \sum_{i=1}^{n}(\tilde{p}_i(x)f_1(L_i^{-1}(x)) + \tilde{q}_i(x)f_2(L_i^{-1}(x)) + \tilde{r}_i(L_i^{-1}(x))) \chi_{I_i}(x), \ x \in I.$$

In virtue of ZIFS, we can construct $2^n$ different fractal interpolation functions by using the same data set and the same vertical scaling factors.

**Remark.** If $\varepsilon = (0, 0, \cdots, 0)$, then the ZHVFIF coinsides with the HVFIF introduced in [6], and if $\beta_i(x) = 0, i \in \mathbf{N}_n$ and all vertical scaling factors are constants, then it coincides with the ZFIF in [9].

The following example shows that the ZHVFIF offers the better approximation than the ZFIF when we interpolate the data set taken from Weierstrass function, the well-known fractal function.

**Example 1.** Let us construct FIFs which interpolates the data set

$$\{(-1, -2), (-0.6667, 0.5), (-0.3333, -0.5), (0, 2), (0.3333, -0.5), (0.6667, 0.5), (1, -2)\},$$

which are taken from Weierstrass function $W(x) = \sum_{s=0}^{\infty} 0.5^s \cos(3^s \pi x), \ x \in [-1, 1]$.

We first construct the affine ZFIF $f$ with the following vertical scaling factors and signature:

$$p = (-0.8, -0.6, -0.7, -0.5, -0.7, -0.6), \ \varepsilon = (0, 1, 0, 1, 0, 1).$$

Next, we construct the ZHVFIF $f_1$ with the extended data set $\{(-1, -2, -4), (-0.6667, 0.5, 1),$ $(-0.3333, -0.5, -1)\}, (0, 2, 3), (0.3333, -0.5, -1.2), (0.6667, 0.5, 0.9), (1, -2, -3)\}$, the vertical scaling

factors

$$\mathbf{p} = (-0.8,\ -0.6,\ -0.7,\ -0.5,\ -0.7,\ -0.6),\quad \mathbf{q} = (0.4,\ 0.3,\ 0.4,\ 0.3,\ 0.3,\ 0.3),$$
$$\widetilde{\mathbf{p}} = (0.8,\ 0.6,\ 0.7,\ 0.6,\ 0.7,\ 0.6),\quad \widetilde{\mathbf{q}} = (-0.4,\ -0.3,\ -0.4,\ -0.3,\ -0.4,\ -0.3)$$

and the same signature $\boldsymbol{\varepsilon}$. Fig. 1 shows the graph of Weierstrass function, $f$ and $f_1$, respectively. We estimate the errors between $f$, $f_1$ and Weierstrass function $W$ using $\dfrac{1}{x_n - x_0}\displaystyle\int_{x_0}^{x_n} |f(x) - W(x)|\,dx$. The errors are 0.979 for $f$ and 0.384 for $f_1$.

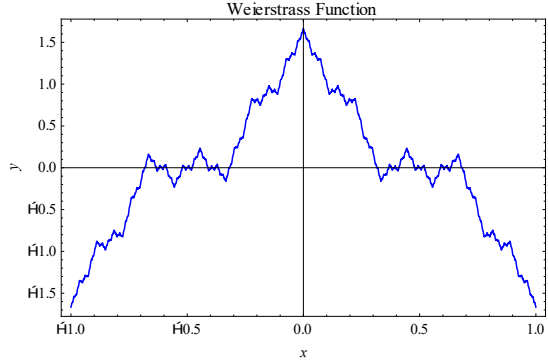

(1) Weierstrass Function

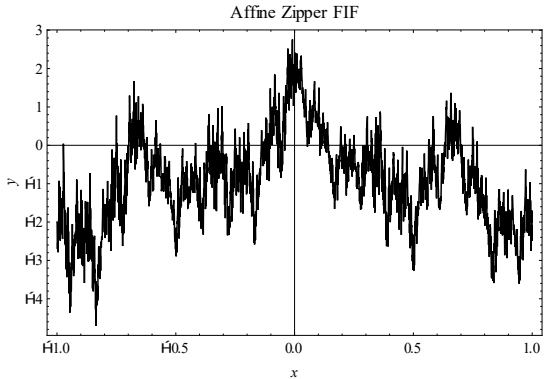

(2) ZFIF

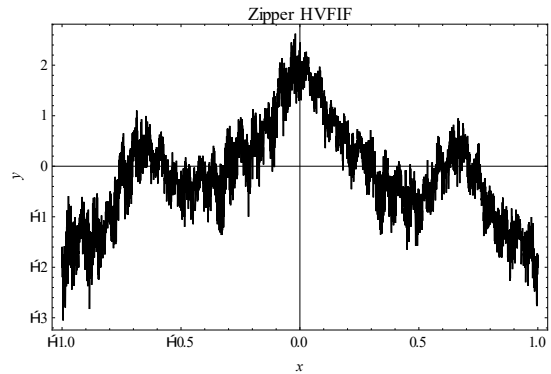

(3) ZHVFIF

Fig. 1 Weierstrass function, ZFIF and ZHVFIF

This example shows that the ZHVFIF gives the better approximation to $W$ than the ZFIF.

## 3. Boundedness of Zipper Hidden Variable Fractal Interpolation Functions

Many real data such as image data, meteorological data, speech data, and electrocardiographic data require that the graphs of interpolation functions are laid inside a desired rectangle. In this section, we find a condition on vertical scaling factors under which the graphs of ZHVFIFs are contained inside the prescribed rectangle containing the data set $P$.

Let $k_1$, $k_2$, $\widetilde{k}_1$ and $\widetilde{k}_2$ be arbitrary real numbers satisfying $k_1 < \min\limits_{i \in \mathbf{N}_n^0}\{y_i\}$, $k_2 > \max\limits_{i \in \mathbf{N}_n^0}\{y_i\}$,

$\tilde{k}_1 < \min\limits_{i \in \mathbf{N}_n^0}\{z_i\}$, $\tilde{k}_2 > \max\limits_{i \in \mathbf{N}_n^0}\{z_i\}$. We choose any $\omega$ satisfying $\overline{S} < \omega < 1$, and denote the maximum value and the minimum value of a function $g$ by $g^{\max}$ and $g_{\min}$, respectively.

**Theorem 3.** *If the maximum and minimum of vertical scaling factors satisfy the following conditions:*

$$p_i^{\min} > -\min\left\{\frac{\min\{y_{i-1},\ y_i\} - k_1}{k_2 - \min\{y_0,\ y_n\}},\ \frac{k_2 - \max\{y_{i-1},\ y_i\}}{\max\{y_0,\ y_n\} - k_1}\right\},$$

$$p_i^{\max} < \min\left\{\frac{\min\{y_{i-1},\ y_i\} - k_1}{\max\{y_0,\ y_n\} - k_1},\ \frac{k_2 - \max\{y_{i-1},\ y_i\}}{k_2 - \min\{y_0,\ y_n\}}\right\},$$

$$q_i^{\min} \geq -\min\left\{\frac{(\min\{y_{i-1},\ y_i\} - k_1) - \max\{p_i^{\max}(y_0 - k_1),\ p_i^{\max}(y_n - k_1),\ -p_i^{\min}(k_2 - y_0),\ -p_i^{\min}(k_2 - y_n)\}}{\tilde{k}_2 - \max\{z_0,\ z_n\}},\right.$$

$$\left.\frac{(k_2 - \max\{y_{i-1},\ y_i\}) - \max\{p_i^{\max}(k_2 - y_0),\ p_i^{\max}(k_2 - y_n),\ -p_i^{\min}(y_0 - k_1),\ -p_i^{\min}(y_n - k_1)\}}{\max\{z_0,\ z_n\} - \tilde{k}_1}\right\},$$

$$q_i^{\max} \leq \min\left\{\frac{(\min\{y_{i-1},\ y_i\} - k_1) - \max\{p_i^{\max}(y_0 - k_1),\ p_i^{\max}(y_n - k_1),\ -p_i^{\min}(k_2 - y_0),\ -p_i^{\min}(k_2 - y_n)\}}{\max\{z_0,\ z_n\} - \tilde{k}_1},\right.$$

$$\left.\frac{(k_2 - \max\{y_{i-1},\ y_i\}) - \max\{p_i^{\max}(k_2 - y_0),\ p_i^{\max}(k_2 - y_n),\ -p_i^{\min}(y_0 - k_1),\ -p_i^{\min}(y_n - k_1)\}}{\tilde{k}_2 - \max\{z_0,\ z_n\}}\right\},$$

$$\tilde{p}_i^{\min} > -\min\left\{\omega - p_i^{\max},\ \omega + p_i^{\min},\ \frac{\min\{y_{i-1},\ y_i\} - k_1}{k_2 - \min\{y_0,\ y_n\}},\ \frac{k_2 - \max\{y_{i-1},\ y_i\}}{\max\{y_0,\ y_n\} - k_1}\right\},$$

$$\tilde{p}_i^{\max} < \min\left\{\omega - p_i^{\max},\ \omega + p_i^{\min},\ \frac{\min\{y_{i-1},\ y_i\} - k_1}{\max\{y_0,\ y_n\} - k_1},\ \frac{k_2 - \max\{y_{i-1},\ y_i\}}{k_2 - \min\{y_0,\ y_n\}}\right\},$$

$$\tilde{q}_i^{\min} \geq -\min\left\{\omega - q_i^{\max},\ \omega + q_i^{\min},\right.$$

$$\frac{(\min\{y_{i-1},\ y_i\} - k_1) - \max\{\tilde{p}_i^{\max}(y_0 - k_1),\ \tilde{p}_i^{\max}(y_n - k_1),\ -\tilde{p}_i^{\min}(k_2 - y_0),\ -\tilde{p}_i^{\min}(k_2 - y_n)\}}{\tilde{k}_2 - \max\{z_0,\ z_n\}},$$

$$\left.\frac{(k_2 - \max\{y_{i-1},\ y_i\}) - \max\{\tilde{p}_i^{\max}(k_2 - y_0),\ \tilde{p}_i^{\max}(k_2 - y_n),\ -\tilde{p}_i^{\min}(y_0 - k_1),\ -\tilde{p}_i^{\min}(y_n - k_1)\}}{\max\{z_0,\ z_n\} - \tilde{k}_1}\right\},$$

$$\widetilde{q}_i^{\max} \leq \min\{\omega - q_i^{\max}, \ \omega + q_i^{\min},$$

$$\frac{(\min\{y_{i-1}, \ y_i\} - k_1) - \max\{\widetilde{p}_i^{\max}(y_0 - k_1), \ \widetilde{p}_i^{\max}(y_n - k_1), \ -\widetilde{p}_i^{\min}(k_2 - y_0), \ -\widetilde{p}_i^{\min}(k_2 - y_n)\}}{\max\{z_0, \ z_n\} - \widetilde{k}_1},$$

$$\left.\frac{(k_2 - \max\{y_{i-1}, \ y_i\}) - \max\{\widetilde{p}_i^{\max}(k_2 - y_0), \ \widetilde{p}_i^{\max}(k_2 - y_n), \ -\widetilde{p}_i^{\min}(y_0 - k_1), \ -\widetilde{p}_i^{\min}(y_n - k_1)\}}{\widetilde{k}_2 - \max\{z_0, \ z_n\}}\right\},$$

*then the graph of ZHVFIF is laid inside the rectangle* $I \times [k_1, \ k_2]$, *i.e.*, $Gr(f_1) \subset I \times [k_1, \ k_2]$.

**Proof.** Firstly, we find a sufficient condition for $p_i$, $q_i$ under which the graph of ZHVFIF $f_1$ is contained in the rectangle $I \times [k_1, \ k_2]$. In other words, we need to find the condition for $F_{i,1}(I \times [k_1, \ k_2] \times [\widetilde{k}_1, \ \widetilde{k}_2]) \subseteq [k_1, \ k_2]$, that is, for any $(x, y, z) \in I \times [k_1, \ k_2] \times [\widetilde{k}_1, \ \widetilde{k}_2]$,

$$k_1 \leq F_{i,1}(x, \ y, \ z) \leq k_2. \tag{2}$$

We divide the proof into 4 parts as follows:

**Case 1.** $p_i(x) \geq 0$, $q_i(x) \geq 0$, $x \in I_i$. In this case, we have

$$\begin{aligned}
F_{i,1}(x, y, z) &= p_i(L_i(x))\{y - (\mu y_n + (1-\mu)y_0)\} + q_i(L_i(x))\{z - (\mu z_n + (1-\mu)z_0)\} \\
&\quad + (\mu y_{i-\varepsilon_i} + (1-\mu)y_{i-1+\varepsilon_i}) \\
&\geq p_i(L_i(x))(k_1 - \max\{y_0, \ y_n\}) + q_i(L_i(x))(\widetilde{k}_1 - \max\{z_0, \ z_n\}) + \min\{y_{i-1}, \ y_i\} \\
&\geq p_i^{\max}(k_1 - \max\{y_0, \ y_n\}) + q_i^{\max}(\widetilde{k}_1 - \max\{z_0, \ z_n\}) + \min\{y_{i-1}, \ y_i\}.
\end{aligned} \tag{3}$$

If we take $p_i^{\max}$ such that

$$p_i^{\max} \leq \frac{\min\{y_{i-1}, \ y_i\} - k_1}{\max\{y_0, \ y_n\} - k_1},$$

then we get

$$\min\{y_{i-1}, \ y_i\} - k_1 - p_i^{\max}(\max\{y_0, \ y_n\} - k_1) \geq 0.$$

Therefore, if for the above $p_i^{\max}$, we take $q_i^{\max}$ satisfying the inequality

$$q_i^{\max} \leq \frac{(\min\{y_{i-1}, \ y_i\} - k_1) - p_i^{\max}(\max\{y_0, \ y_n\} - k_1)}{\max\{z_0, \ z_n\} - \widetilde{k}_1},$$

then we obtain

$$p_i^{\max}(k_1 - \max\{y_0, \ y_n\}) + \min\{y_{i-1}, \ y_i\} - k_1 \geq q_i^{\max}(\max\{z_0, \ z_n\} - \widetilde{k}_1) \geq 0.$$

Therefore, from inequality (3), we obtain $F_{i,1}(x, y, z) \geq k_1$.

Similarly, we have

$$F_{i,1}(x, y, z) \leq p_i(L_i(x))(k_2 - \min\{y_0, y_n\}) + q_i(L_i(x))(\tilde{k}_2 - \min\{z_0, z_n\}) + \max\{y_{i-1}, y_i\}$$
$$\leq p_i^{\max}(k_2 - \min\{y_0, y_n\}) + q_i^{\max}(\tilde{k}_2 - \min\{z_0, z_n\}) + \max\{y_{i-1}, y_i\}. \quad (4)$$

If we take $p_i^{\max}$ such that

$$p_i^{\max} \leq \frac{k_2 - \max\{y_{i-1}, y_i\}}{k_2 - \min\{y_0, y_n\}},$$

then we obtain

$$k_2 - \max\{y_{i-1}, y_i\} - p_i^{\max}(k_2 - \min\{y_0, y_n\}) \geq 0.$$

Therefore, if for the above $p_i^{\max}$, we take $q_i^{\max}$ satisfying the inequality

$$q_i^{\max} \leq \frac{(k_2 - \max\{y_{i-1}, y_i\}) - p_i^{\max}(k_2 - \min\{y_0, y_n\})}{\tilde{k}_2 - \min\{z_0, z_n\}},$$

then we get

$$k_2 - \max\{y_{i-1}, y_i\} - p_i^{\max}(k_2 - \min\{y_0, y_n\}) \geq q_i^{\max}(\tilde{k}_2 - \min\{z_0, z_n\}) \geq 0.$$

Therefore, from inequality (4), we get $F_{i,1}(x, y, z) \leq k_2$.

Hence, the condition for the inequality (2) is as follows:

$$p_i^{\max} \leq \min\left\{ \frac{\min\{y_{i-1}, y_i\} - k_1}{\max\{y_0, y_n\} - k_1}, \frac{k_2 - \max\{y_{i-1}, y_i\}}{k_2 - \min\{y_0, y_n\}} \right\},$$

$$q_i^{\max} \leq \min\left\{ \frac{(\min\{y_{i-1}, y_i\} - k_1) - p_i^{\max}(\max\{y_0, y_n\} - k_1)}{\max\{z_0, z_n\} - \tilde{k}_1}, \right.$$

$$\left. \frac{(k_2 - \max\{y_{i-1}, y_i\}) - p_i^{\max}(k_2 - \min\{y_0, y_n\})}{\tilde{k}_2 - \min\{z_0, z_n\}} \right\}.$$

**Case 2.** $p_i(x) < 0, q_i(x) < 0, x \in I_i$. In this case, as in Case 1, we have

$$F_{i,1}(x, y, z) \geq p_i^{\min}(k_2 - \min\{y_0, y_n\}) + q_i^{\min}(\tilde{k}_2 - \min\{z_0, z_n\}) + \min\{y_{i-1}, y_i\} \quad (5)$$

and

$$F_{i,1}(x, y, z) \leq -p_i^{\min}(\max\{y_0, y_n\} - k_1) - q_i^{\min}(\max\{z_0, z_n\} - \widetilde{k}_1) + \max\{y_{i-1}, y_i\}. \tag{6}$$

If we take $p_i^{\min}$ and $q_i^{\min}$ such that

$$p_i^{\min} > -\frac{\min\{y_{i-1}, y_i\} - k_1}{k_2 - \min\{y_0, y_n\}}, \quad q_i^{\min} \geq -\frac{(\min\{y_{i-1}, y_i\} - k_1) + p_i^{\min}(k_2 - \min\{y_0, y_n\})}{\widetilde{k}_2 - \min\{z_0, z_n\}},$$

then we obtain

$$k_1 - \min\{y_{i-1}, y_i\} - p_i^{\min}(k_2 - \min\{y_0, y_n\}) \leq q_i^{\min}(\widetilde{k}_2 - \min\{z_0, z_n\}) < 0.$$

Therefore, from inequality (5), we get $F_{i,1}(x, y, z) \geq k_1$.

Similarly, if we take $p_i^{\min}$ and $q_i^{\min}$ such that

$$p_i^{\min} > -\frac{k_2 - \max\{y_{i-1}, y_i\}}{\max\{y_0, y_n\} - k_1}, \quad q_i^{\min} \geq -\frac{(k_2 - \max\{y_{i-1}, y_i\}) + p_i^{\min}(\max\{y_0, y_n\} - k_1)}{\max\{z_0, z_n\} - \widetilde{k}_1},$$

then we get

$$\max\{y_{i-1}, y_i\} - k_2 - p_i^{\min}(\max\{y_0, y_n\} - k_1) \leq q_i^{\min}(\max\{z_0, z_n\} - \widetilde{k}_1) < 0.$$

Furthermore, from inequality (6), we get $F_{i,1}(x, y, z) \leq k_2$.

Hence, the condition for the inequality (2) is as follows:

$$p_i^{\min} > -\min\left\{\frac{\min\{y_{i-1}, y_i\} - k_1}{k_2 - \min\{y_0, y_n\}}, \frac{k_2 - \max\{y_{i-1}, y_i\}}{\max\{y_0, y_n\} - k_1}\right\},$$

$$q_i^{\min} \geq -\min\left\{\frac{(\min\{y_{i-1}, y_i\} - k_1) + p_i^{\min}(k_2 - \min\{y_0, y_n\})}{\widetilde{k}_2 - \min\{z_0, z_n\}},\right.$$

$$\left.\frac{(k_2 - \max\{y_{i-1}, y_i\}) + p_i^{\min}(\max\{y_0, y_n\} - k_1)}{\max\{z_0, z_n\} - \widetilde{k}_1}\right\}.$$

**Case 3.** $p_i(x) \geq 0$, $q_i(x) < 0$, $x \in I_i$. In this case, we have

$$F_{i,1}(x, y, z) \geq -p_i^{\max}(\max\{y_0, y_n\} - k_1) + q_i^{\min}(\widetilde{k}_2 - \min\{z_0, z_n\}) + \min\{y_{i-1}, y_i\},$$

$$F_{i,1}(x, y, z) \leq p_i^{\max}(k_2 - \min\{y_0, y_n\}) - q_i^{\max}(\max\{z_0, z_n\} - \widetilde{k}_1) + \max\{y_{i-1}, y_i\}.$$

As in Case 1, if we take $p_i^{\max}$ and $q_i^{\min}$ such that

$$p_i^{\max} < \min\left\{\frac{\min\{y_{i-1},\ y_i\} - k_1}{\max\{y_0,\ y_n\} - k_1},\ \frac{k_2 - \max\{y_{i-1},\ y_i\}}{k_2 - \min\{y_0,\ y_n\}}\right\}$$

and

$$q_i^{\min} \geq -\min\left\{\frac{(\min\{y_{i-1},\ y_i\} - k_1) - p_i^{\max}(\max\{y_0,\ y_n\} - k_1)}{\widetilde{k}_2 - \min\{z_0,\ z_n\}},\right.$$

$$\left.\frac{(k_2 - \max\{y_{i-1},\ y_i\}) - p_i^{\max}(k_2 - \min\{y_0,\ y_n\})}{\max\{z_0,\ z_n\} - \widetilde{k}_1}\right\},$$

then we get $k_1 \leq F_{i,1}(x, y, z) \leq k_2$.

**Case 4.** $p_i(x) < 0,\ q_i(x) \geq 0,\ x \in I_i$. In this case, we have

$$F_{i,1}(x, y, z) \geq p_i^{\min}(k_2 - \min\{y_0,\ y_n\}) - q_i^{\max}(\max\{z_0,\ z_n\} - \widetilde{k}_1) + \min\{y_{i-1},\ y_i\}$$

and

$$F_{i,1}(x, y, z) \leq -p_i^{\min}(\max\{y_0,\ y_n\} - k_1) + q_i^{\max}(\widetilde{k}_2 - \min\{z_0,\ z_n\}) + \max\{y_{i-1},\ y_i\}.$$

Furthermore, as in Case 1, if we take

$$p_i^{\min} \geq -\min\left\{\frac{\min\{y_{i-1},\ y_i\} - k_1}{k_2 - \min\{y_0,\ y_n\}},\ \frac{k_2 - \max\{y_{i-1},\ y_i\}}{\max\{y_0,\ y_n\} - k_1}\right\}$$

and

$$q_i^{\max} \leq \min\left\{\frac{(\min\{y_{i-1},\ y_i\} - k_1) + p_i^{\min}(k_2 - \min\{y_0,\ y_n\})}{\max\{z_0,\ z_n\} - \widetilde{k}_1},\right.$$

$$\left.\frac{(k_2 - \max\{y_{i-1},\ y_i\}) + p_i^{\min}(\max\{y_0,\ y_n\} - k_1)}{\widetilde{k}_2 - \min\{z_0,\ z_n\}}\right\},$$

then we get $k_1 \leq F_{i,1}(x, y, z) \leq k_2$.

Hence, if $p_i, q_i, i \in \mathbf{N}_n$ satisfy the conditions of the theorem, then $Gr(f_1) \subseteq I \times [k_1, k_2]$.

Next, in the similar way, we can obtain the sufficient condition for $\widetilde{p}_i,\ \widetilde{q}_i$ under which the graph of ZFIF $f_2$ is contained in the rectangle $I \times [\widetilde{k}_1,\ \widetilde{k}_2]$. Furthermore under the condition $\overline{S} < \omega < 1$, we can obtain the assumptions in the theorem. □

**Example 2.** The data set is given as follows:

$$\{(-1,-4),(-0.6667,2),(-0.3333,0),(0,3),(0.3333,1),(0.6667,4),(1,-2)\}.$$

Let us extend the data set as follows:

$$P=\{(-1,-4,-3),(-0.6667,2,1),(-0.3333,0,-2),(0,3,4),(0.3333,1,3),(0.6667,4,6),(1,-2,2)\}.$$

Now, we construct the ZHVFIF interpolating $P_0$ whose graph is in $[-1,1]\times[-10,10]$.

Let $\varepsilon=(0,1,0,1,0,1)$. According to Theorem 3, if $k_1=-10, k_2=10, \widetilde{k}_1=-6$ and $\widetilde{k}_2=8$, then vertical scaling factors must satisfy the following conditions:

$$-0.428\leq p_1\leq 0.571,\ -0.714\leq p_2\leq 0.571,\ -0.714\leq p_3\leq 0.5,$$
$$-0.785\leq p_4\leq 0.5,\ -0.75\leq p_5\leq 0.428,\ -0.571\leq p_6\leq 0.428,$$
$$-0.067\leq q_1\leq 0.05,\ -0.125\leq q_2\leq 0.333,\ -0.033\leq q_3\leq 0.025,$$
$$-0.088\leq q_4\leq 0.117,\ -0.15\leq q_5\leq 0.2,\ -0.05\leq q_6\leq 0.07,$$
$$-0.428\leq \widetilde{p}_1\leq 0.5,\ -0.4\leq \widetilde{p}_2\leq 0.4,\ -0.2\leq \widetilde{p}_3\leq 0.2,$$
$$-0.45\leq \widetilde{p}_4\leq 0.45,\ -0.3\leq \widetilde{p}_5\leq 0.3,\ -0.5\leq \widetilde{p}_6\leq 0.428,$$
$$-0.667\leq \widetilde{q}_1\leq 0.05,\ -0.125\leq \widetilde{q}_2\leq 0.333,\ -0.033\leq \widetilde{q}_3\leq 0.025,$$
$$-0.088\leq q_4\leq 0.117,\ -0.15\leq q_5\leq 0.2,\ -0.05\leq q_6\leq 0.067.$$

On the basis of these conditions, we take the following vertical scaling factors and construct the ZHVFIF:

$$\mathbf{p}=(0.5,\ -0.18,\ 0.53,\ -0.17,\ 0.52,\ 0.48),$$
$$\mathbf{q}=(0.17,\ -0.054,\ 0.053,\ -0.034,\ 0.23,\ 0.14),$$
$$\widetilde{\mathbf{p}}=(0.36,\ -0.21,\ 0.33,\ -0.18,\ 0.34,\ 0.37),$$
$$\widetilde{\mathbf{q}}=(0.29,\ -0.05,\ 0.22,\ -0.14,\ 0.38,\ 0.23).$$

The graph of the ZHVFIF is laid inside $[-1,1]\times[-10,10]$, which is shown in the Fig. 2.

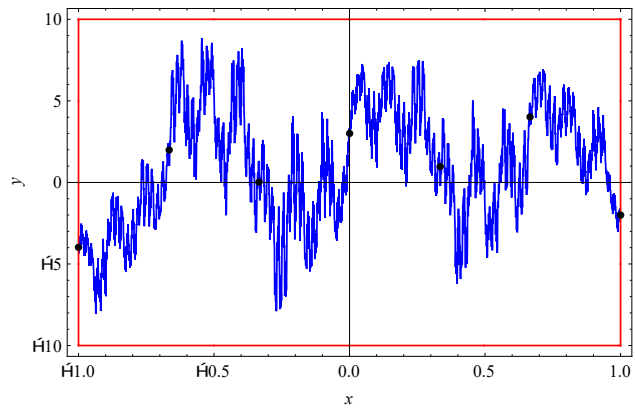

Fig. 2 Bounded ZHVFIF with a signature $\varepsilon=(0,1,0,1,0,1)$.

## 4. Positivity of Zipper Hidden Variable Fractal Interpolation Functions

There are many real data whose values are always positive such as wind speed or relative humidity. In this section, we consider a sufficient condition for $p_i$, $q_i$ under which the ZHVFIF $f_1$ has positivity.

Suppose that $y_i \geq 0$, $z_i \geq 0$, $i \in \mathbf{N}_n$ in the extended data set $P$. Let us denote

$$\omega_i := \max\{\|p_i\|_\infty, \|q_i\|_\infty\}, \quad \tilde{\omega}_i := \max\{\|\tilde{p}_i\|_\infty, \|\tilde{q}_i\|_\infty\},$$

$$M := \left(\max_{i \in \mathbf{N}_n^0}\{y_i\} + \max_{i \in \mathbf{N}_n^0}\{z_i\} - \max_{i \in \mathbf{N}_n}\{\omega_i + \tilde{\omega}_i\} \cdot (\min\{y_0, y_n\} + \min\{z_0, z_n\})\right) \Big/ \left(1 - \max_{i \in \mathbf{N}_n}\{\omega_i + \tilde{\omega}_i\}\right).$$

Then, we can easily see that for the constructed ZHVFIF $f_1$ and ZFIF $f_2$, we have $|f_1(x)| \leq M$, $|f_2(x)| \leq M$, $x \in I$.

**Theorem 4.** *If the maximum and minimum of vertical scaling factors satisfy the following conditions:*

$$p_i^{\min} > -\frac{\min\{y_{i-1}, y_i\}}{M - \min\{y_0, y_n\}}, \quad p_i^{\max} < \frac{\min\{y_{i-1}, y_i\}}{\max\{y_0, y_n\}},$$

$$q_i^{\min} \geq -\min\left\{\frac{\min\{y_{i-1}, y_i\} + p_i^{\min}(M - \min\{y_0, y_n\})}{M - \min\{z_0, z_n\}}, \frac{\min\{y_{i-1}, y_i\} - p_i^{\max}\min\{y_0, y_n\}}{M - \min\{z_0, z_n\}}\right\},$$

$$q_i^{\max} \leq \min\left\{\frac{\min\{y_{i-1}, y_i\} - p_i^{\max}\max\{y_0, y_n\}}{\max\{z_0, z_n\}}, \frac{\min\{y_{i-1}, y_i\} + p_i^{\min}(M - \max\{y_0, y_n\})}{\max\{z_0, z_n\}}\right\},$$

$$\tilde{p}_i^{\min} > -\min\left\{\omega - p_i^{\max}, \omega + p_i^{\min}, \frac{\min\{y_{i-1}, y_i\}}{M - \min\{y_0, y_n\}}\right\},$$

$$\tilde{p}_i^{\max} < \min\left\{\omega - p_i^{\max}, \omega + p_i^{\min}, \frac{\min\{y_{i-1}, y_i\}}{\max\{y_0, y_n\}}\right\},$$

$$\tilde{q}_i^{\min} \geq -\min\left\{\omega - q_i^{\max}, \omega + q_i^{\min}, \frac{\min\{y_{i-1}, y_i\} + \tilde{p}_i^{\min}(M - \min\{y_0, y_n\})}{M - \min\{z_0, z_n\}}, \frac{\min\{y_{i-1}, y_i\} - \tilde{p}_i^{\max}\min\{y_0, y_n\}}{M - \min\{z_0, z_n\}}\right\},$$

$$\tilde{q}_i^{\max} \leq \min\left\{\omega - q_i^{\max}, \omega + q_i^{\min}, \frac{\min\{y_{i-1}, y_i\} - \tilde{p}_i^{\max}\max\{y_0, y_n\}}{\max\{z_0, z_n\}}, \frac{\min\{y_{i-1}, y_i\} + \tilde{p}_i^{\min}(M - \max\{y_0, y_n\})}{\max\{z_0, z_n\}}\right\},$$

*then the ZHVFIF $f_1$ and the ZFIF $f_2$ are positive, i.e. $f_1(x) \geq 0$, $f_2(x) \geq 0$, $x \in I$.*

**Proof.** Firstly, we have to find the sufficient conditions under which for $x \in I$, $y \geq 0$, $z \geq 0$, we get

$$F_{i,1}(x, y, z) \geq 0. \tag{7}$$

**Case 1.** $p_i(x) \geq 0$, $q_i(x) \geq 0$, $x \in I_i$. In this case, we have

$$F_{i,1}(x, y, z) = p_i(L_i(x))\{y - (\mu y_n + (1-\mu)y_0)\} + q_i(L_i(x))\{z - (\mu z_n + (1-\mu)z_0)\}$$
$$+ (\mu y_{i-\varepsilon_i} + (1-\mu)y_{i-1+\varepsilon_i}) \quad (8)$$
$$\geq p_i(L_i(x))(y - \max\{y_0, y_n\}) + q_i(L_i(x))(z - \max\{z_0, z_n\}) + \min\{y_{i-1}, y_i\}$$
$$\geq -p_i^{\max} \cdot \max\{y_0, y_n\}) - q_i^{\max} \cdot \max\{z_0, z_n\}) + \min\{y_{i-1}, y_i\}.$$

If we take $p_i^{\max}$ such that $p_i^{\max} \leq \dfrac{\min\{y_{i-1}, y_i\}}{\max\{y_0, y_n\}}$, then we get

$$\min\{y_{i-1}, y_i\} - p_i^{\max} \max\{y_0, y_n\} \geq 0.$$

Therefore, if for the above $p_i^{\max}$, we take $q_i^{\max}$ satisfying the inequality

$$q_i^{\max} \leq \dfrac{\min\{y_{i-1}, y_i\} - p_i^{\max} \max\{y_0, y_n\}}{\max\{z_0, z_n\}},$$

then we get

$$-q_i^{\max} \cdot \max\{z_0, z_n\}) + \min\{y_{i-1}, y_i\} \geq p_i^{\max} \cdot \max\{y_0, y_n\}) \geq 0.$$

Therefore, from inequality (8), we get $F_{i,1}(x, y, z) \geq 0$. Hence, the condition for the inequality (7) is as follows:

$$p_i^{\max} \leq \dfrac{\min\{y_{i-1}, y_i\}}{\max\{y_0, y_n\}}, \quad q_i^{\max} \leq \dfrac{\min\{y_{i-1}, y_i\} - p_i^{\max} \max\{y_0, y_n\}}{\max\{z_0, z_n\}}.$$

**Case 2.** $p_i(x) < 0$, $q_i(x) < 0$, $x \in I_i$. In this case, as in Case 1, we have

$$F_{i,1}(x, y, z) \geq p_i(L_i(x))(M - \min\{y_0, y_n\}) + q_i(L_i(x))(M - \min\{z_0, z_n\}) + \min\{y_{i-1}, y_i\}.$$

If we take $p_i^{\min}$ and $q_i^{\min}$ such that

$$p_i^{\min} \geq -\dfrac{\min\{y_{i-1}, y_i\}}{M - \min\{y_0, y_n\}}, \quad q_i^{\min} \geq -\dfrac{\min\{y_{i-1}, y_i\} + p_i^{\min}(M - \min\{y_0, y_n\})}{M - \min\{z_0, z_n\}},$$

then we obtain

$$F_{i,1}(x, y, z) \geq p_i(L_i(x))(M - \min\{y_0, y_n\}) + q_i(L_i(x))(M - \min\{z_0, z_n\}) + \min\{y_{i-1}, y_i\} \geq 0.$$

**Case 3.** $p_i(x) \geq 0$, $q_i(x) < 0$, $x \in I_i$. In this case, we have

$$F_{i,1}(x, y, z) \geq -p_i(L_i(x))\max\{y_0, y_n\} + q_i(L_i(x))(M - \min\{z_0, z_n\}) + \min\{y_{i-1}, y_i\}.$$

If we take $p_i^{\min}$ and $q_i^{\min}$ such that

$$p_i^{\max} \leq \frac{\min\{y_{i-1}, y_i\}}{\max\{y_0, y_n\}}, \quad q_i^{\min} \geq -\frac{\min\{y_{i-1}, y_i\} - p_i^{\max}\min\{y_0, y_n\}}{M - \min\{z_0, z_n\}},$$

then from the above inequality, we obtain $F_{i,1}(x, y, z) \geq 0$.

**Case 4.** $p_i(x) < 0$, $q_i(x) \geq 0$, $x \in I_i$. In this case, we have

$$F_{i,1}(x, y, z) \geq p_i(L_i(x))(M - \min\{y_0, y_n\}) - q_i(L_i(x))\max\{z_0, z_n\} + \min\{y_{i-1}, y_i\}.$$

If we take $p_i^{\min}$ and $q_i^{\max}$ such that

$$p_i^{\min} \geq -\frac{\min\{y_{i-1}, y_i\}}{M - \min\{y_0, y_n\}}, \quad q_i^{\max} \leq \frac{\min\{y_{i-1}, y_i\} + p_i^{\min}(M - \max\{y_0, y_n\})}{\max\{z_0, z_n\}},$$

then from the above inequality, inequality (7) holds.

Hence, we can see that $f_1(x) \geq 0$, $x \in I$ under the conditions of the theorem on vertical scaling factors $p_i, q_i$.

Similarly, we can find the sufficient conditions such that $f_2(x) \geq 0$, $x \in I$. □

**Example 3.** The data set is given as follows:
$$\{(-1, 2), (-0.6667, 1.4), (-0.3333, 1.5), (0, 1.2), (0.3333, 2.1), (0.6667, 1.6), (1, 1.3)\}.$$

Now we extend the data set as follows:
$$P = \{(-1, 2, 1), (-0.6667, 1.4, 1.6), (-0.3333, 1.5, 4), (0, 1.2, 0.3),$$
$$(0.3333, 2.1, 2), (0.6667, 1.6, 4.2), (1, 1.3, 2.1)\}.$$

We take a signature $\varepsilon = (0, 1, 0, 1, 0, 1)$. By Theorem 4, for the ZHVFIF to be positive, the vertical scaling factors $p_i$, $i \in \mathbf{N}_6$ must satisfy the following conditions:
$$-0.245 \leq p_1 \leq 0.7, \quad -0.245 \leq p_2 \leq 0.7, \quad -0.21 \leq p_3 \leq 0.6,$$
$$-0.21 \leq p_4 \leq 0.6, \quad -0.28 \leq p_5 \leq 0.8, \quad -0.228 \leq p_6 \leq 0.65.$$

If we take $p_1 = 0.5$, then Theorem 4 yields
$$-0.125 \leq q_1 \leq 0.19, \quad -0.245 \leq \tilde{p}_1 \leq 0.4, \quad -0.155 \leq \tilde{q}_1 \leq 0.324.$$

Therefore, we take $q_1 = 0.17$, $\tilde{p}_1 = 0.36$, $\tilde{q}_1 = 0.29$ under the above condition.

Similarly, we take the following:

$p_2 = -0.18$, $q_2 = -0.054$, $\tilde{p}_2 = -0.21$, $\tilde{q}_2 = -0.05$ such that
$$-0.062 \leq q_2 \leq 0.238, \quad -0.245 \leq \tilde{p}_2 \leq 0.4, \quad -0.062 \leq \tilde{q}_2 \leq 0.238,$$

$p_3 = 0.53$, $q_3 = 0.053$, $\tilde{p}_3 = 0.33$, $\tilde{q}_3 = 0.22$ such that
$$-0.085 \leq q_3 \leq 0.067, \quad -0.21 \leq \tilde{p}_3 \leq 0.37, \quad -0.128 \leq \tilde{q}_3 \leq 0.257,$$

$p_4 = -0.17$, $q_4 = -0.034$, $\tilde{p}_4 = -0.18$, $\tilde{q}_4 = -0.14$ such that
$$-0.038 \leq q_4 \leq 0.166, \quad -0.211 \leq \tilde{p}_4 \leq 0.6, \quad -0.161 \leq \tilde{q}_4 \leq 0.143,$$

$p_5 = 0.52$, $q_5 = 0.018$, $\widetilde{p}_5 = 0.76$, $\widetilde{q}_5 = 0.88$ such that
$$-0.154 \leq q_5 \leq 0.267, \quad -0.28 \leq \widetilde{p}_5 \leq 0.38, \quad -0.193 \leq \widetilde{q}_5 \leq 0.438,$$
$p_6 = 0.48$, $q_6 = 0.14$, $\widetilde{p}_6 = 0.37$, $\widetilde{q}_6 = 0.23$ such that
$$-0.112 \leq q_5 \leq 0.161, \quad -0.228 \leq \widetilde{p}_5 \leq 0.42, \quad -0.198 \leq \widetilde{q}_5 \leq 0.266.$$

Hence, we get the following vertical scaling factors:

$\mathbf{p} = (0.5, -0.18, 0.53, -0.17, 0.52, 0.48)$, $\mathbf{q} = (0.17, -0.054, 0.053, -0.034, 0.23, 0.14)$,

$\widetilde{\mathbf{p}} = (0.36, -0.21, 0.33, -0.18, 0.34, 0.37)$, $\widetilde{\mathbf{q}} = (0.29, -0.05, 0.22, -0.14, 0.38, 0.23)$

and construct the ZHVFIF using them. Then, the ZHVFIF has positivity and Fig. 3 shows the graph of the ZHVFIF.

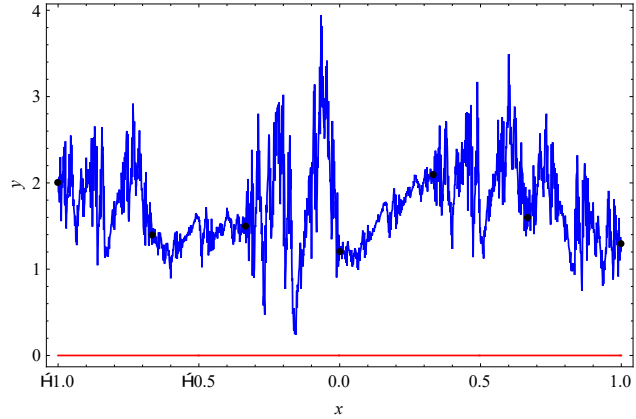

Fig. 3 Positive ZHVFIF with a signature $\varepsilon = (0, 1, 0, 1, 0, 1)$.

## 5. Piecewise Slop Preservation of Zipper Hidden Variable Fractal Interpolation Functions

Interpolation of many real data such as meteorology data, earthquake data or time series in physiology requires the fact that the graph of the interpolation function is laid between pairs of parallel straight lines, which is called piecewise slop preservation. In this section, we find the conditions on vertical scaling factors for the ZHVFIFs to preserve piecewise slop.

Let us denote

$$y_{\min} := \min_{i \in \mathbf{N}_n}\{m_i x_{i-1} + b_{i1}, \; m_n x_n + b_{n1}\}, \quad y_{\max} := \max_{i \in \mathbf{N}_n}\{m_i x_{i-1} + b_{i1}, \; m_n x_n + b_{n1}\}, \quad m_i := \frac{y_i - y_{i-1}}{x_i - x_{i-1}}.$$

**1) The case where a data set is given either above or below the piecewise lines**

Let the piecewise lines be given by $y = m_i x + b_{i1}$, $i \in \mathbf{N}_n$ such that a data set is laid above the piecewise lines, i.e., $m_i x + b_{i1} < y_i$, $m_1 x_0 + b_{11} < y_0$. We find the conditions on vertical scaling factors under which the graph of the constructed ZHVFIF is laid above the piecewise lines, i.e. for any $i \in \mathbf{N}_n$, $m_i x + b_{i1} \leq f_1(x), \; x \in I_i$.

Let us denote

$$\widetilde{k} := \frac{\max_{i \in \mathbf{N}_0^n}\{|y_i + z_i|\} - \omega \min\{|y_0 + z_0|, |y_n + z_n|\} - (1-\omega)y_{\min}}{1-\omega}.$$

Then, we can easily see that $|f_2(x)| \le \widetilde{k}$, $x \in I$.

**Theorem 5.** *Let the maximum and minimum of vertical scaling factors satisfy the following conditions:*

$$p_i^{\min} > -\frac{y_0 - m_1 x_0 - b_{11}}{U - \min\{y_0, y_n\}}, \quad p_i^{\max} < \frac{y_0 - m_1 x_0 - b_{11}}{\max\{y_0, y_n\} - y_{\min}},$$

$$q_i^{\min} \ge -\min\left\{\frac{(y_0 - m_1 x_0 - b_{11}) + p_i^{\min}(U - \min\{y_0, y_n\})}{\widetilde{k} - \min\{z_0, z_n\}}, \frac{(y_0 - m_1 x_0 - b_{11}) - p_i^{\max}(\max\{y_0, y_n\} - y_{\min})}{\widetilde{k} - \min\{z_0, z_n\}}\right\},$$

$$q_i^{\max} \le \min\left\{\frac{(y_0 - m_1 x_0 - b_{11}) - p_i^{\max}(\max\{y_0, y_n\} - y_{\min})}{\max\{z_0, z_n\} + \widetilde{k}}, \frac{(y_0 - m_1 x_0 - b_{11}) + p_i^{\min}(M - \min\{y_0, y_n\})}{\max\{z_0, z_n\} + \widetilde{k}}\right\},$$

$$\widetilde{p}_i^{\min} \ge -\min\{\omega + p_i^{\min}, \omega - p_i^{\max}\}, \quad \widetilde{p}_i^{\max} \le \min\{\omega + p_i^{\min}, \omega - p_i^{\max}\},$$

$$\widetilde{q}_i^{\min} \ge -\min\{\omega + q_i^{\min}, \omega - q_i^{\max}\}, \quad \widetilde{q}_i^{\max} \le \min\{\omega + q_i^{\min}, \omega - q_i^{\max}\}.$$

*Then the graph of the ZHVFIF $f_1$ is constrained above piecewise lines $y = m_i x + b_{i1}$, $i \in \mathbf{N}_n$.*

**Proof.** Firstly, we find a sufficient condition such that for any $(x, y, z) \in I \times [y_{\min}, U] \times [-\widetilde{k}, \widetilde{k}]$,

$$F_{i,1}(x, y, z) \ge m_i L_i(x) + b_{i1}. \tag{9}$$

**Case 1.** $p_i(x) \ge 0$, $q_i(x) \ge 0$, $x \in I_i$. In this case, we have

$$F_{i,1}(x, y, z) = p_i(L_i(x))\{y - [\mu y_n + (1-\mu)y_0]\} + q_i(L_i(x))\{z - [\mu z_n + (1-\mu)z_0]\} +$$

$$+ [\mu y_{i-\varepsilon_i} + (1-\mu)y_{i-1+\varepsilon_i}]$$

$$\ge -p_i(L_i(x))(\max\{y_0, y_n\} - y_{\min}) - q_i(L_i(x))(\max\{z_0, z_n\} + \widetilde{k}) + [\mu y_{i-\varepsilon_i} + (1-\mu)y_{i-1+\varepsilon_i}]$$

$$\ge -p_i^{\max}(\max\{y_0, y_n\} - y_{\min}) - q_i^{\max}(\max\{z_0, z_n\} + \widetilde{k}) + [\mu y_{i-\varepsilon_i} + (1-\mu)y_{i-1+\varepsilon_i}]. \tag{10}$$

We know that

$$[\mu y_{i-\varepsilon_i} + (1-\mu)y_{i-1+\varepsilon_i}] - [m_i L_i(x) + b_{i1}] = y_{i-1+\varepsilon_i} - m_i x_{i-1+\varepsilon_i} - b_{i1} = y_0 - m_1 x_0 - b_{11}.$$

From this, if we take $p_i^{\max}$ such that

$$p_i^{\max} < \frac{y_0 - m_1 x_0 - b_{11}}{\max\{y_0, y_n\} - y_{\min}},$$

then we get

$$(y_0 - m_1 x_0 - b_{11}) - p_i^{\max}(\max\{y_0, y_n\} - y_{\min}) \geq 0.$$

Therefore, if for the above $p_i^{\max}$, we take $q_i^{\max}$ that satisfies the inequality

$$q_i^{\max} \leq \frac{(y_0 - m_1 x_0 - b_{11}) - p_i^{\max}(\max\{y_0, y_n\} - y_{\min})}{\max\{z_0, z_n\} + \widetilde{k}},$$

then we obtain

$$0 \leq q_i^{\max}(\max\{z_0, z_n\} + \widetilde{k}) \leq -p_i^{\max}(\max\{y_0, y_n\} - y_{\min}) + y_0 - m_1 x_0 - b_{11}.$$

Therefore, from inequality (10), we get $F_{i,1}(x, y, z) \geq m_i L_i(x) + b_{i1}$.

Hence, the condition for the inequality (9) is as follows:

$$p_i^{\max} < \frac{y_0 - m_1 x_0 - b_{11}}{\max\{y_0, y_n\} - y_{\min}}, \quad q_i^{\max} \leq \frac{(y_0 - m_1 x_0 - b_{11}) - p_i^{\max}(\max\{y_0, y_n\} - y_{\min})}{\max\{z_0, z_n\} + \widetilde{k}}$$

**Case 2.** $p_i(x) < 0$, $q_i(x) < 0$, $x \in I_i$. We can consider this case in the similar way in Case 1. If we have

$$p_i^{\min} > -\frac{y_0 - m_1 x_0 - b_{11}}{U - \min\{y_0, y_n\}}, \quad q_i^{\min} \geq -\frac{(y_0 - m_1 x_0 - b_{11}) + p_i^{\min}(U - \min\{y_0, y_n\})}{\widetilde{k} - \min\{z_0, z_n\}},$$

then inequality (9) holds.

**Case 3.** $p_i(x) \geq 0$, $q_i(x) < 0$, $x \in I_i$. If we take $p_i^{\max}, q_i^{\min}$ such that

$$p_i^{\max} \leq \frac{y_0 - m_1 x_0 - b_{11}}{\max\{y_0, y_n\} - y_{\min}}, \quad q_i^{\min} \geq -\frac{(y_0 - m_1 x_0 - b_{11}) - p_i^{\max}(\max\{y_0, y_n\} - y_{\min})}{\widetilde{k} - \min\{z_0, z_n\}},$$

then the inequality (9) holds.

**Case 4.** $p_i(x) < 0$, $q_i(x) \geq 0$, $x \in I_i$. In this case, the conditions for the inequality (9) is as follows:

$$p_i^{\min} > -\frac{y_0 - m_1 x_0 - b_{11}}{U - \min\{y_0, y_n\}}, \quad q_i^{\max} \leq \frac{(y_0 - m_1 x_0 - b_{11}) + p_i^{\min}(M - \min\{y_0, y_n\})}{\max\{z_0, z_n\} - \widetilde{k}}.$$

Thus, if $p_i^{\min}$, $p_i^{\max}$, $q_i^{\min}$, $q_i^{\max}$ $i \in \mathbf{N}_n$ satisfy the assumptions in the theorem, then in any subinterval $I_i$, the graph of ZHVFIF $f_1$ is laid above a straight line $y = m_i x + b_{i1}$.

In the similar way, we can get the conditions for $\widetilde{p}_i^{\min}$, $\widetilde{p}_i^{\max}$, $\widetilde{q}_i^{\min}$, $\widetilde{q}_i^{\max}$ under the condition $\overline{S} < \omega < 1$. □

**Example 4.** The data set is given as follows:
$$\{(-1, -4), (-0.6667, 2), (-0.3333, 0), (0, 3), (0.3333, 4), (0.6667, 4), (1, -2)\}.$$
We extend it as follows:
$$P = \{(-1, -4, -3), (-0.667, 2, 1), (-0.333, 0, -2), (0, 3, 4), (0.3333, 1, 3), (0.6667, 4, 6), (1, -2, 2)\}.$$
We take piecewise lines below the data set as follows:
$$y = \begin{cases} 18x + 12, & -1 \leq x < -0.667 \\ -6x - 4, & -0.667 \leq x < -0.333 \\ 9x + 1, & -0.333 \leq x < 0 \\ -6x + 1, & 0 \leq x < 0.333 \\ 9x - 4, & 0.333 \leq x < 0.667 \\ -18x + 14, & 0.667 \leq x \leq 1 \end{cases}$$

Let us choose a signature and vertical scaling factors to satisfy Theorem 5 as follows:
$$\varepsilon = (0, 1, 0, 1, 0, 1),$$
$$\mathbf{p} = (-0.1, 0.48, -0.09, 0.45, -0.08, 0.4),$$
$$\mathbf{q} = (-0.00096, 0.0007, -0.002, 0.0017, -0.0045, 0.0035),$$
$$\tilde{\mathbf{p}} = (-0.78, 0.41, -0.79, 0.44, -0.81, 0.48),$$
$$\tilde{\mathbf{q}} = (-0.899, 0.8991, -0.895, 0.897, 0.893, 0.895).$$

Fig. 4 shows the graph of the constructed ZHVFIF.

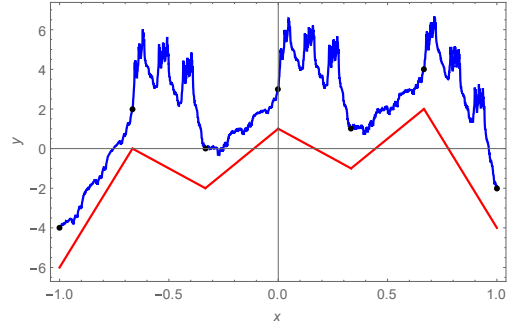

Fig. 4 ZHVFIF constrained below

Similarly, we can find the conditions under which for any $i \in \mathbf{N}_n$, we have $f_1(x) \leq m_i x + b_{i2}$, $x \in I_i$ in the case where a data set is laid below the piecewise lines $y = m_i x + b_{i2}, x \in I_i$ such that $y_i < m_i x + b_{i2}, x \in I_i$, $y_0 < b_{12} + m_1 x_0$. Let us denote
$$\tilde{\ell} := \frac{\max_{i \in \mathbf{N}_0^n}\{|y_i + z_i|\} - \omega \min\{|y_0 + z_0|, |y_n + z_n|\} - (1-\omega)U}{1-\omega}.$$

Then, in this case, we can easily get $|f_2(x)| \leq \tilde{\ell}, x \in I$.

**Theorem 6.** *If the maximum and minimum of vertical scaling factors satisfy the following conditions, then the graph of ZHVFIF $f_1$ is constrained below piecewise lines $y = m_i x + b_{i2}, i \in \mathbf{N}_n$.*

$$p_i^{\min} > -\frac{m_1 x_0 + b_{12} - y_0}{U + \max\{y_0, y_n\}}, \quad p_i^{\max} < \frac{m_1 x_0 + b_{12} - y_0}{y_{\max} - \min\{y_0, y_n\}},$$

$$q_i^{\min} = -\min\left\{\frac{(m_1 x_0 + b_{12} - y_0) + p_i^{\min}(U + \max\{y_0, y_n\})}{\max\{z_0, z_n\} + \widetilde{\ell}}, \frac{(m_1 x_0 + b_{12} - y_0) - p_i^{\max}(y_{\max} - \min\{y_0, y_n\})}{\max\{z_0, z_n\} + \widetilde{\ell}}\right\},$$

$$q_i^{\max} = \min\left\{\frac{(m_1 x_0 + b_{12} - y_0) + p_i^{\min}(\max\{y_0, y_n\} + U)}{\widetilde{\ell} - \min\{z_0, z_n\}}, \frac{(m_1 x_0 + b_{12} - y_0) - p_i^{\max}(y_{\max} - \min\{y_0, y_n\})}{\widetilde{\ell} - \min\{z_0, z_n\}}\right\},$$

$$\widetilde{p}_i^{\min} = -\min\{\omega + p_i^{\min}, \omega - p_i^{\max}\}, \quad \widetilde{p}_i^{\max} = \min\{\omega + p_i^{\min}, \omega - p_i^{\max}\},$$

$$\widetilde{q}_i^{\min} = -\min\{\omega + q_i^{\min}, \omega - q_i^{\max}\}, \quad \widetilde{q}_i^{\max} = \min\{\omega + q_i^{\min}, \omega - q_i^{\max}\}.$$

**Example 5.** The data set and the extended data set are given as in example 4. We take the piecewise lines above the data set as follows:

$$y = \begin{cases} 18x + 16, & -1 \le x < -0.667 \\ -6x, & -0.667 \le x < -0.333 \\ 9x + 5, & -0.333 \le x < 0 \\ -6x + 5, & 0 \le x < 0.333 \\ 9x, & 0.333 \le x < 0.667 \\ -18x + 18, & 0.667 \le x \le 1 \end{cases}$$

Let us choose a signature and vertical scaling factors to satisfy Theorem 6 as follows:

$$\boldsymbol{\varepsilon} = (1, 0, 1, 0, 1, 0),$$

$\mathbf{p} = (0.1, -0.48, 0.09, -0.45, 0.08, -0.4)$, $\mathbf{q} = (0.00096, -0.0007, 0.002, -0.0017, 0.0045, -0.0035)$,
$\widetilde{\mathbf{p}} = (0.78, -0.41, 0.79, -0.44\ 0.81, -0.48)$, $\widetilde{\mathbf{q}} = (0.899, -0.8991, 0.895, -0.897, 0.893, -0.895)$.

The graph of the constructed ZHVFIF is below the piecewise lines, which is shown in Figure 5.

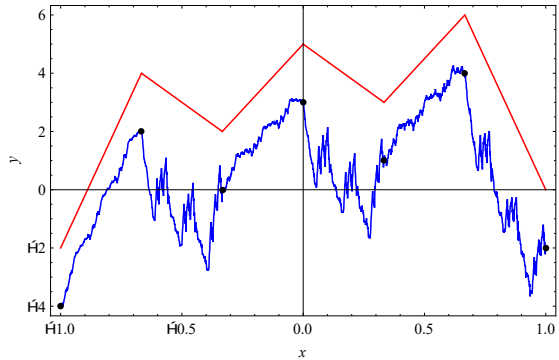

Fig. 5 ZHVFIF constrained above

**2) The case where piecewise lines below and above a data set are given**

Let the piecewise lines be given by $y = m_i x + b_{i1}$, $y = m_i x + b_{i2}$, $x \in I_i$, $i \in \mathbf{N}_n$ such that a data set is contained between $y = m_i x + b_{i1}$ and $y = m_i x + b_{i2}$, $x \in I_i$, i.e.

$$m_i x + b_{i1} \le y_i \le m_i x + b_{i2}, \quad x \in I_i, \quad m_1 x_0 + b_{11} < y_0 < m_1 x_0 + b_{12}.$$

Now, we find the conditions under which for any $i \in \mathbf{N}_n$, $m_i x + b_{i1} \le f_1(x) \le m_i x + b_{i2}$, $x \in I_i$. In this case, we can easily see that $|f_2(x)| \leq \widetilde{k}$, $x \in I$.

**Theorem 7.** *If the maximum and minimum of vertical scaling factors satisfy the following conditions, then the graph of the ZHVFIF $f_1$ is contained between the piecewise lines $y = m_i x + b_{i1}$ and $y = m_i x + b_{i2}$:*

$$p_i^{\min} > -\min\left\{\frac{y_0 - m_i x_0 - b_{11}}{y_{\max} - \min\{y_0, y_n\}}, \frac{m_1 x_0 + b_{12} - y_0}{\max\{y_0, y_n\} - y_{\min}}\right\},$$

$$p_i^{\max} < \min\left\{\frac{y_0 - m_i x_0 - b_{11}}{\max\{y_0, y_n\} - y_{\min}}, \frac{m_1 x_0 + b_{12} - y_0}{y_{\max} - \min\{y_0, y_n\}}\right\},$$

$$q_i^{\min} = -\min\left\{\frac{(y_0 - m_i x_0 - b_{11}) + p_i^{\min}(y_{\max} - \min\{y_0, y_n\})}{\widetilde{k} - \min\{z_0, z_n\}}, \frac{(m_1 x_0 + b_{12} - y_0) + p_i^{\min}(\max\{y_0, y_n\} - y_{\min})}{\max\{z_0, z_n\} - \widetilde{k}}\right.$$

$$\left.\frac{(y_0 - m_i x_0 - b_{11}) - p_i^{\max}(\max\{y_0, y_n\} - y_{\min})}{\widetilde{k} - \min\{z_0, z_n\}}, \frac{(m_1 x_0 + b_{12} - y_0) - p_i^{\max}(y_{\max} - \min\{y_0, y_n\})}{\max\{z_0, z_n\} - \widetilde{k}}\right\},$$

$$q_i^{\max} = \min\left\{\frac{(y_0 - m_i x_0 - b_{11}) - p_i^{\max}(\max\{y_0, y_n\} - y_{\min})}{\max\{z_0, z_n\} - \widetilde{k}}, \frac{(m_1 x_0 + b_{12} - y_0) - p_i^{\max}(y_{\max} - \min\{y_0, y_n\})}{\widetilde{k} - \min\{z_0, z_n\}}\right.$$

$$\left.\frac{(y_0 - m_i x_0 - b_{11}) + p_i^{\min}(y_{\max} - \min\{y_0, y_n\})}{\max\{z_0, z_n\} - \widetilde{k}}, \frac{(m_1 x_0 + b_{12} - y_0) + p_i^{\min}(\max\{y_0, y_n\} - y_{\min})}{\widetilde{k} - \min\{z_0, z_n\}}\right\},$$

$$\widetilde{p}_i^{\min} = -\min\{\omega + p_i^{\min}, \omega - p_i^{\max}\}, \quad \widetilde{p}_i^{\max} = \min\{\omega + p_i^{\min}, \omega - p_i^{\max}\},$$

$$\widetilde{q}_i^{\min} = -\min\{\omega + q_i^{\min}, \omega - q_i^{\max}\}, \quad \widetilde{q}_i^{\max} = \min\{\omega + q_i^{\min}, \omega - q_i^{\max}\}.$$

**Proof.** Firstly, we find the conditions under which for any $x \in I$, $y \in [y_{\min}, y_{\max}]$ and $z \in [-\widetilde{k}, \widetilde{k}]$, we get

$$m_i L_i(x) + b_{i1} \leq F_{i,1}(x, y, z) \leq m_i L_i(x) + b_{i2}. \tag{11}$$

**Case 1.** $p_i(x) \geq 0$, $q_i(x) \geq 0$, $x \in I_i$. In this case, we have

$$F_{i,1}(x, y, z) = p_i(L_i(x))\{y - [\mu y_n + (1-\mu)y_0]\} + q_i(L_i(x))\{z - [\mu z_n + (1-\mu)z_0]\} +$$
$$+ [\mu y_{i-\varepsilon_i} + (1-\mu) y_{i-1+\varepsilon_i}] \tag{12}$$
$$\geq -p_i^{\max}(\max\{y_0, y_n\} - y_{\min}) - q_i^{\max}(\max\{z_0, z_n\} + \widetilde{k}) + [\mu y_{i-\varepsilon_i} + (1-\mu) y_{i-1+\varepsilon_i}].$$

We know that

$$[\mu y_{i-\varepsilon_i} + (1-\mu)y_{i-1+\varepsilon_i}] - [m_i L_i(x) + b_{i1}] = y_{i-1+\varepsilon_i} - m_i x_{i-1+\varepsilon_i} - b_{i1} = y_0 - m_i x_0 - b_{11},$$

$$[m_i L_i(x) + b_{i2}] - [\mu y_{i-\varepsilon_i} + (1-\mu)y_{i-1+\varepsilon_i}] = m_1 x_0 + b_{12} - y_0.$$

If we take $p_i^{\max}$ such that

$$p_i^{\max} \leq \frac{y_0 - m_i x_0 - b_{11}}{\max\{y_0,\ y_n\} - y_{\min}},$$

then we can easily see $(y_0 - m_i x_0 - b_{11}) - p_i^{\max}(\max\{y_0,\ y_n\} - y_{\min}) \geq 0$. Therefore, if for the above $p_i^{\max}$, we take $q_i^{\max}$ satisfying the inequality

$$q_i^{\max} \leq \frac{(y_0 - m_i x_0 - b_{11}) - p_i^{\max}(\max\{y_0,\ y_n\} - y_{\min})}{\max\{z_0,\ z_n\} + \widetilde{k}},$$

then we obtain

$$0 \leq q_i^{\max}(\max\{z_0,\ z_n\} + \widetilde{k}) \leq -p_i^{\max}(\max\{y_0,\ y_n\} - y_{\min}) + y_0 - m_i x_0 - b_{11}.$$

Therefore, from inequality (12), we obtain $F_{i,1}(x, y, z) \geq m_i L_i(x) + b_{i1}$.

Similarly, we have

$$F_{i,1}(x, y, z) \leq p_i^{\max}(y_{\max} - \min\{y_0,\ y_n\}) + q_i^{\max}(\widetilde{k} - \min\{z_0,\ z_n\}) + [\mu y_{i-\varepsilon_i} + (1-\mu)y_{i-1+\varepsilon_i}]. \quad (13)$$

If we take $p_i^{\max}$ and $q_i^{\max}$ such that

$$p_i^{\max} \leq \frac{m_1 x_0 + b_{12} - y_0}{y_{\max} - \min\{y_0,\ y_n\}},\quad q_i^{\max} \leq \frac{(m_1 x_0 + b_{12} - y_0) - p_i^{\max}(y_{\max} - \min\{y_0,\ y_n\})}{\widetilde{k} - \min\{z_0,\ z_n\}},$$

then we obtain

$$m_1 x_0 + b_{12} - y_0 - p_i^{\max}(y_{\max} - \min\{y_0,\ y_n\}) \geq q_i^{\max}(\widetilde{k} - \min\{z_0,\ z_n\}) \geq 0.$$

Therefore, from inequality (13), we get $F_{i,1}(x, y, z) \leq m_i L_i(x) + b_{i2}$.

Hence, if the following conditions hold, then (11) holds:

$$p_i^{\max} \leq \min\left\{\frac{y_0 - m_i x_0 - b_{11}}{\max\{y_0,\ y_n\} - y_{\min}},\ \frac{m_1 x_0 + b_{12} - y_0}{y_{\max} - \min\{y_0,\ y_n\}}\right\},$$

$$q_i^{\max} \leq \min\left\{\frac{(y_0 - m_i x_0 - b_{11}) - p_i^{\max}(\max\{y_0, y_n\} - y_{\min})}{\max\{z_0, z_n\} + \widetilde{k}}, \frac{(m_1 x_0 + b_{12} - y_0) - p_i^{\max}(y_{\max} - \min\{y_0, y_n\})}{\widetilde{k} - \min\{z_0, z_n\}}\right\}.$$

**Case 2.** $p_i(x) < 0$, $q_i(x) < 0$, $x \in I_i$. We can deal with this case in the similar way in Case 1.

If we take $p_i^{\min}$ and $q_i^{\min}$ such that

$$p_i^{\min} \geq -\min\left\{\frac{y_0 - m_i x_0 - b_{11}}{y_{\max} - \min\{y_0, y_n\}}, \frac{m_1 x_0 + b_{12} - y_0}{\max\{y_0, y_n\} - y_{\min}}\right\},$$

$$q_i^{\min} \geq -\min\left\{\frac{(y_0 - m_i x_0 - b_{11}) + p_i^{\min}(y_{\max} - \min\{y_0, y_n\})}{\widetilde{k} - \min\{z_0, z_n\}}, \frac{(m_1 x_0 + b_{12} - y_0) + p_i^{\min}(\max\{y_0, y_n\} - y_{\min})}{\max\{z_0, z_n\} + \widetilde{k}}\right\},$$

then (11) holds.

**Case 3.** $p_i(x) \geq 0$, $q_i(x) < 0$, $x \in I_i$. If we take $p_i^{\min}$ and $q_i^{\min}$ such that

$$p_i^{\max} \leq \min\left\{\frac{y_0 - m_i x_0 - b_{11}}{\max\{y_0, y_n\} - y_{\min}}, \frac{m_1 x_0 + b_{12} - y_0}{y_{\max} - \min\{y_0, y_n\}}\right\},$$

$$q_i^{\min} \geq -\min\left\{\frac{(y_0 - m_i x_0 - b_{11}) - p_i^{\max}(\max\{y_0, y_n\} - y_{\min})}{\widetilde{k} - \min\{z_0, z_n\}}, \frac{(m_1 x_0 + b_{12} - y_0) - p_i^{\max}(y_{\max} - \min\{y_0, y_n\})}{\max\{z_0, z_n\} + \widetilde{k}}\right\},$$

then (11) holds.

**Case 4.** $p_i(x) < 0$, $q_i(x) \geq 0$, $x \in I_i$. In this case, the condition for the inequality (11) is as follows:

$$p_i^{\min} \geq -\min\left\{\frac{y_0 - m_i x_0 - b_{11}}{y_{\max} - \min\{y_0, y_n\}}, \frac{m_1 x_0 + b_{12} - y_0}{\max\{y_0, y_n\} - y_{\min}}\right\},$$

$$q_i^{\max} \leq \min\left\{\frac{(y_0 - m_i x_0 - b_{11}) + p_i^{\min}(y_{\max} - \min\{y_0, y_n\})}{\max\{z_0, z_n\} + \widetilde{k}}, \frac{(m_1 x_0 + b_{12} - y_0) + p_i^{\min}(\max\{y_0, y_n\} - y_{\min})}{\widetilde{k} - \min\{z_0, z_n\}}\right\}.$$

Thus, if $p_i$, $q_i$, $i \in \mathbf{N}_n$ satisfy the assumptions in the theorem, then in any subinterval $I_i$, the graph of the ZHVFIF $f_1$ lies between piecewise lines $y = m_i x + b_{i1}$ and $y = m_i x + b_{i2}$, $x \in I_i$.

Next, in the similar way, we can get the conditions for $\widetilde{p}_i^{\min}$, $\widetilde{p}_i^{\max}$, $\widetilde{q}_i^{\min}$, $\widetilde{q}_i^{\max}$ under the condition $\overline{S} < \omega < 1$. □

**Example 6.** The data set is given as follows:

$\{(-1, 5), (-0.6667, -2), (-0.3333, 3), (0, 0), (0.3333, 1), (0.6667, -3), (1, 2)\}$.

We extend it as follows:
$$P = \{(-1, 5, 2), (-0.6667, -2, -4), (-0.3333, 3, -2), (0, 0, -3),$$
$$(0.3333, 1, 2), (0.6667, -3, 0), (1, 2, 5)\}.$$

We take a pair of piecewise lines above and below the data set as follows:

$$y = \begin{cases} -21x - 17.5, & -1 \leq x < -0.667 \\ 15x + 6.5, & -0.667 \leq x < -0.333 \\ -9x - 1.5, & -0.333 \leq x < 0 \\ 3x - 1.5, & 0 \leq x < 0.333 \\ -12x + 3.5, & 0.333 \leq x < 0.667 \\ 15x - 14.5, & 0.667 \leq x \leq 1 \end{cases}, \quad y = \begin{cases} -21x - 14.5, & -1 \leq x < -0.667 \\ 15x + 9.5, & -0.667 \leq x < -0.333 \\ -9x + 1.5, & -0.333 \leq x < 0 \\ 3x + 1.5, & 0 \leq x < 0.333 \\ -12x + 6.5, & 0.333 \leq x < 0.667 \\ 15x - 11.5, & 0.667 \leq x \leq 1 \end{cases}.$$

Take $\varepsilon = (0, 1, 0, 1, 0, 1)$. By Theorem 7, the vertical scaling factors $p_i$, $i \in \mathbf{N}_6$ have to satisfy the following condition:
$$-0.158 \leq p_i \leq 0.158, \; i \in \mathbf{N}_6.$$

From this condition, first, we take $p_1 = 0.15$. Then from Theorem 7, we have
$$-0.0187 \leq q_1 \leq 0.0068, \; -0.75 \leq \tilde{p}_1 \leq 0.75, \; -0.894 \leq \tilde{q}_1 \leq 0.894,$$
and therefore we take $q_1 = 0.006, \tilde{p}_1 = 0.7, \tilde{q}_1 = 0.89$.

Similarly, we take

$p_2 = -0.14$ and $q_2 = -0.015$, $\tilde{p}_2 = -0.75$, $\tilde{q}_2 = -0.85$ such that
$$-0.1545 \leq q_2 \leq 0.0425, \; -0.76 \leq \tilde{p}_2 \leq 0.76, \; -0.885 \leq \tilde{q}_2 \leq 0.885,$$

$p_3 = 0.1355$ and $q_3 = 0.018$, $\tilde{p}_3 = 0.74$, $\tilde{q}_3 = 0.8$ such that
$$-0.081 \leq q_3 \leq 0.02, \; -0.7645 \leq \tilde{p}_3 \leq 0.7645, \; -0.882 \leq \tilde{q}_3 \leq 0.882,$$

$p_4 = -0.13$ and $q_4 = -0.02$, $\tilde{p}_4 = -0.7$, $\tilde{q}_4 = -0.8$ such that
$$-0.024 \leq q_4 \leq 0.083, \; -0.77 \leq \tilde{p}_4 \leq 0.77, \; -0.88 \leq \tilde{q}_4 \leq 0.88,$$

$p_5 = 0.135$ and $q_5 = 0.018$, $\tilde{p}_5 = 0.76$, $\tilde{q}_5 = 0.88$ such that
$$-0.054 \leq q_5 \leq 0.0197, \; -0.765 \leq \tilde{p}_5 \leq 0.765, \; -0.882 \leq \tilde{q}_5 \leq 0.882,$$

$p_6 = -0.155$ and $q_6 = -0.002$, $\tilde{p}_6 = -0.74$, $\tilde{q}_6 = -0.89$ such that
$$-0.0025 \leq q_6 \leq 0.007, \; -0.745 \leq \tilde{p}_6 \leq 0.745, \; -0.898 \leq \tilde{q}_6 \leq 0.898.$$

Hence, we have the following vertical scaling factors:
$$\mathbf{p} = (0.15, -0.14, 0.1355, -0.13, 0.135, -0.155),$$
$$\mathbf{q} = (0.006, -0.015, 0.018, -0.02, 0.018, -0.002),$$
$$\tilde{\mathbf{p}} = (0.7, -0.75, 0.74, -0.7, 0.76, -0.74),$$
$$\tilde{\mathbf{q}} = (0.89, -0.85, 0.8, -0.8, 0.88, -0.89).$$

Then, the graph of the constructed ZHVFIF is laid between two piecewise lines, which is shown in Fig. 6.

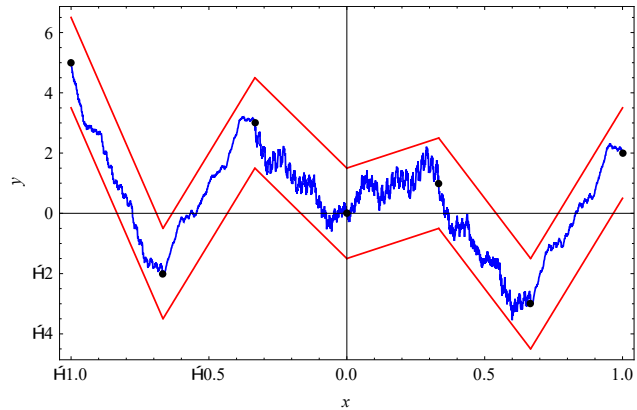

Fig. 6 Piecewise slop preserving ZHVFIF with a signature $\varepsilon = (0, 1, 0, 1, 0, 1)$.

## Conclusion

In this paper, we constructed zipper hidden variable fractal interpolation functions with the zipper. By the virtue of zipper, we can get more flexible and approximative fractal interpolation functions with the same data set. Next, we found some sufficient conditions on the vertical scaling factors for the zipper hidden variable fractal interpolation function to preserve boundedness or positivity, piecewise slop of the data set, which are very important to model natural objects by fractal interpolation functions. In the future, we are going to study fractional calculus of zipper fractal interpolation functions.

**Declaration of Competing Interest**

The authors declare that they have no known competing financial interests or personal relationships that could have appeared to influence the work reported in this paper.

**Data availability**

The data and code that support the findings of this study are available upon reasonable request from the corresponding author.

**CRediT authorship contribution statement**

**Chol Hui Yun:** Conceptualization, Methodology, Project administration, Supervision, Writing-original draft.
**Yu Jong Pak:** Formal analysis, Investigation, Methodology, Software, Writing-review & editing.
**Mi Gyong Ri:** Investigation, Validation, Visualization, Software, Writing-review & editing
**Gyong Ju Ri:** Investigation, Software, Validation, Writing-review & editing

**Funding**

Not applicable.

# References

[1] V.V. Assev, A.V. Tetenov, A.S. Kravchenko, On self-similar Jordan curves on the plane. Sib. Math. J., 44(3) (2003) 379–386.

[2] N. Attia, M. Balegh, R. Amami, R. Amami, On the Fractal interpolation functions associated with Matkowski contractions, Electronic Research Archive, 31(8) (2023) 4652–4668.

[3] B. Balázs, K. Gergely, K. István, Pointwise regularity of parameterized affine zipper fractal curves, Nonlinearity, 31 (2018) 1705–1733.

[4] M.F. Barnsley, Fractal functions and interpolation, Constr. Approx., 2 (1986) 303–329.

[5] M.F. Barnsley, J.H. Elton, Recurrent iterated function systems, Constr. Approx., 5 (1989) 3-31.

[6] M.F. Barnsley, D. Hardin, P. Massopust, Hidden variable fractal interpolation functions, SIAM J. Math. Anal., 20 (1989) 1218-1242.

[7] P. Bouboulis, L. Dalla, Hidden variable vector valued fractal interpolation functions, Fractals, 13(3) (2005) 227–232.

[8] P. Bouboulis, L. Dalla, V. Drakopoulos, Construction of recurrent bivariate fractal interpolation surfaces and computation of their box-counting dimension, J. Approx. Th., 141 (2006), 99-117.

[9] A.K.B. Chand, N.Vijender, P. Viswanathan, A. V. Tetenov, Affine zipper fractal interpolation functions, BIT, 60(2) (2020) 319–344.

[10] G.P. Kapoor, S.A. Prasad, Smoothness of hidden variable bivariate coalescence fractal interpolation surfaces, Int. J. Bifurc, Chaos, 19(7) (2009) 2321–2333.

[11] G.P. Kapoor, S.A. Prasad, Stability of Coalescence Hidden Variable Fractal Interpolation Surfaces, Int. J. Nonlinear Sci., 9(3) (2010) 265-275.

[12] S.K. Katiyar, A.K.B. Chand, G.Saravana Kumar, A new class of rational cubic spline fractal interpolation function and its constrained aspects, Appl.Math.Comput., 346 (2019) 319-335.

[13] J.M. Kim, H.M. Mun, Nonlinear recurrent hidden variable fractal interpolation curves with function vertical scaling factors, Fractals, 28(2020) 2050096.

[14] J.M. Kim, H.J. Kim, H.M. Mun, Nonlinear fractal interpolation curves with function vertical scaling factors. Indian J. Pure Appl. Math., 51(2) (2020) 483-499.

[15] H.J. Kim, J.M. Kim, H.M. Mun, New nonlinear recurrent hidden variable fractal interpolation surfaces, Fractals, 28(2) (2020) 2050038.

[16] W. Metzler, C.H. Yun, Construction of fractal interpolation surfaces on rectangular grids, Int. J. of Bifurcat. Chaos, 20(12) (2010), 4079–4086.

[17] M. Radu, R. Pasupathi, Contractive Multivariate Zipper Fractal Interpolation Functions, Results Math., 79(4) (2024) DOI: 10.1007/s00025-024-02177-5

[18] S. Ri, New types of fractal interpolation surfaces, Chaos Solitons Fractals, 119 (2019) 291–297.

[19] M.G. Ri, C.H. Yun, Smoothness and fractional integral of hidden variable recurrent fractal interpolation function with function vertical scaling factors, Fractals, 29(6) (2021) 2150136, 1-17

[20] M.G. Ri, C.H. Yun, M.H. Kim, Construction of cubic spline hidden variable recurrent fractal interpolation function and its fractional calculus, Chaos Solitons Fractals, 150 (2021) 111177, 1-11.

[21] M.G. Ri, C.H. Yun, Smoothness and fractional integral of hidden variable recurrent fractal interpolation function with function vertical scaling factors, Fractals, 29(6) (2021), 2150136, 1-17.

[22] N.V. Vijay , A.K.B. Chand, C1-Positivity preserving Bi-quintic blended rational quartic zipper fractal interpolation surfaces, Chaos Solitons Fractals, 188 (2024) 115472, DOI: https://doi.org/10.1016/j.chaos.2024.115472.

[23] N.V. Vijay, M.G.P. Prasad, Gurunathan, A novel class of zipper fractal Bezier curves and its graphics applications, Chaos Solitons Fractals, 190 (2024) 115793, DOI: https://doi.org/10.1016/j.chaos.2024.115793.

[24] N.V. Vijay, A.K.B. Chand, Generalized zipper fractal approximation and parameter identification problems, Comp. App. Math., 41 (2022) 155, 1-23.

[25] C.H. Yun, H.C. Choi, H.C. O, Construction of recurrent fractal interpolation surfaces with function vertical scaling factors and estimation of box-counting dimension on rectangular grids, Fractals, 23(3) (2015) 1-10.

[26] C.H. Yun, M.G. Ri, Box-counting dimension and analytic properties of hidden variable fractal interpolation functions with function contractivity factors, Chaos Solitons Fractals, 134 (2020) 109700, 1-10.

[27] C.H. Yun, M.G. Ri, Analytic properties of hidden variable bivariable fractal interpolation functions with four function contractivity factors, Fractals, 27(2) (2019) 1950018, 1-16.